\input amstex.tex
\documentstyle{amsppt}
\magnification=\magstep1
%
%
\def\section#1{\par\bigpagebreak\par\bigpagebreak
    \csname subhead\endcsname #1\endsubhead\par\bigpagebreak}
\def\subsection#1{\par\medpagebreak
    \csname subsubhead\endcsname #1 \endsubsubhead}
\def\Theorem#1#2{\csname proclaim\endcsname{Theorem #1} #2 
    \endproclaim}
\def\Corollary#1{\csname proclaim\endcsname{Corollary} #1 \endproclaim}
\def\Proposition#1#2{\csname proclaim\endcsname{Proposition #1} #2 
    \endproclaim}
\def\Lemma#1#2{\csname proclaim\endcsname{Lemma #1} #2 \endproclaim}
\def\Remark#1#2{\remark{Remark {\rm #1}} #2 \endremark}
\def\Proof#1{\demo{Proof} #1\qed\enddemo}
\def\Proofof#1#2{\demo{Proof of #1} #2\qed\enddemo}
%
%
\define\ep{\epsilon}
\define\dt{\delta}
\define\Dt{\Delta}
\define\ld{\lambda}

\define\Z{{\Bbb Z}}
\define\Q{{\Bbb Q}}

\define\K{{\Bbb K}}
\define\N{{\Bbb N}}
\define\scpr#1#2{\langle{#1},\,{#2}\rangle}

\define\bb#1{\text{\bf #1}}
\define\oT{\overline{T}}
\define\CP{{\Cal P}}
\leftheadtext\nofrills{A.N.\,Kirillov and M.\,Noumi}
\rightheadtext\nofrills{Raising operators for Macdonald polynomials }
\topmatter
\title\nofrills
Affine Hecke algebras and raising operators \\
for Macdonald polynomials 
\endtitle
\author
Anatol N.\, KIRILLOV
\footnote"$^{\ast1}$"
{\hbox{Department of Mathematical Sciences, University of Tokyo;}
\newline
\hbox{\quad St.\,Petersburg Branch of the Steklov Mathematical Institute}}
{\rm and} 
Masatoshi NOUMI
\footnote"$^{\ast2}$"
{\hbox{Department of Mathematics, Kobe University}}
\endauthor
\address
A.N.\,Kirillov: 
\newline\indent
{\eightpoint
Department of Mathematical Sciences, University of Tokyo,
\newline\indent
Komaba, Meguro, Tokyo 153, Japan;
\newline\indent
Steklov Mathematical Institute, 
\newline\indent
Fontanka 27, St.\,Petersburg 191011, Russia
}
\newline\indent
M.\,Noumi: 
\newline\indent
{\eightpoint
Department of Mathematics, Faculty of Science, Kobe University,
\newline\indent
Rokko, Kobe 657, Japan}
\endaddress
\email 
kirillov\@ker.c.u-tokyo.ac.jp\quad
noumi\@math.s.kobe-u.ac.jp
\endemail
\endtopmatter
\document
\section{Introduction}
In this paper we introduce certain raising operators and lowering operators 
for Macdonald polynomials (of type $A_{n-1}$) by means of the Dunkl operators 
due to I.\,Cherednik. 
The raising operators we discuss below are a natural $q$-analogue
of the raising operators for Jack polynomials introduced by 
L.\,Lapointe and L.\, Vinet \cite{LV1, LV2}. 
As an application of our raising operators, we will prove the integrality 
of double Kostka coefficients which had been conjectured by I.G.\,Macdonald 
\cite{Ma1} (apart from the positivity conjecture). 
We will also include some application to a double analogue of the multinomial 
coefficients. 
\par\medpagebreak
Let $\K=\Q(q,t)$ be the field of rational functions in two indeterminates $(q,t)$ 
and $\K[x]^W$ the algebra of symmetric polynomials in 
$n$ variables $x=(x_1,\cdots,x_n)$ over $\K$, $W$ being the symmetric group 
$\frak{S}_n$ of degree $n$. 
The {\it Macdonald polynomials} $P_\ld(x)=P_\ld(x;q,t)$ (or 
{\it symmetric functions with two parameters}, in the terminology of 
I.G.\,Macdonald \cite{Ma1}), 
are a family of symmetric polynomials parametrized by partitions, 
and they form a $\K$-basis of $\K[x]^W$.  
One way to characterize these polynomials is, among others, to consider  
the joint eigenfunctions in $\K[x]^W$ for the commuting family of 
$q$-difference operators
$$
D^{(r)}_x=
t^\binom{r}{2}\sum\Sb I\subset[1,n] \\ |I|=r\endSb
\prod\Sb i\in I\\ j\not\in I\endSb \frac{tx_i-x_j}{x_i-x_j} 
\prod_{i\in I}T_{q,x_i}
\quad (r=0,1,\cdots, n).  
\tag{1}
$$
The Macdonald polynomial $P_\ld(x)$ is characterized as the joint eigenfunction
of $D^{(r)}_x$ ($r=0,1,\cdots,n$) that has the leading term $m_\ld(x)$ 
under the dominance order of partitions when it is expressed as a linear 
combination of monomial symmetric functions $m_\mu(x)$.  
As in \cite{Ma1} (VI.8.3), we also use another normalization 
$J_\ld(x)=c_\ld P_\ld(x)$,
called the ``integral form'' of $P_\ld(x)$. 
\par
We define the operators $B_m$ and $A_m$ ($m=0,1,\cdots,n$), involving 
$q$-shift operators and permutations, as
$$
\align
&B_m=\sum_{k_1<k_2<\cdots<k_m}
x_{k_1}\cdots x_{k_m} (1-t^m Y_{k_1})(1-t^{m-1}Y_{k_2})
\cdots(1-tY_{k_m}),
\tag{2}\\
&
A_m=\sum_{k_1<k_2<\cdots<k_m}
\frac{1}{x_{k_1}\cdots x_{k_m}}(1-Y^\ast_{k_1})(1-tY^\ast_{k_2})
\cdots(1-t^{m-1}Y^\ast_{k_m}),
\endalign
$$ 
by means of the {\it Dunkl operators} $Y_{k}$ and their
dual version  $Y^\ast_{k}$  ($k=1,\cdots,n$).  
The Dunkl operators for Macdonald polynomials, defined in terms of 
a representation of the (extended) affine Hecke algebra, are due to 
Cherednik \cite{C}. 
(For our normalization of $Y_k$ and $Y^\ast_k$, see Section 2.) 
The main result of this paper is:
\Theorem{1}{
The operators $B_m$ and $A_m$ are raising and lowering operators, 
respectively, in the sense that
$$
B_m J_\ld(x)=J_{\ld+(1^m)}(x), \quad\text{and}\quad 
A_mJ_\ld(x)
=a_\ld J_{\ld-(1^m)}(x),
\tag{3}
$$
for each partition $\ld$ with $\ell(\ld)\le m$, 
with some $a_\ld\in\Z[q,t]$.  
}
\noindent(See Theorem 3.1 in Section 3.) 
\par
Theorem 1 implies that, 
for any partition $\ld=(\ld_1,\cdots,\ld_n)$, 
the Macdonald polynomial $J_\ld(x)$ 
is obtained by a successive application of the operators $B_m$ starting 
from $J_{0}(x)=1$:
$$
J_{\ld}(x)=(B_n)^{\ld_n}(B_{n-1})^{\ld_{n-1}-\ld_n}\cdots(B_1)^{\ld_1-\ld_2}(1).  
\tag{4}
$$
We can make use of this expression to study transition coefficients between 
Macdonald polynomials $J_\ld(x)=J_\ld(x;q,t)$
and other symmetric functions. 
In particular we have
\Theorem{2}{
For any partition $\ld$ and $\mu$,  
the double Kostka coefficient $K_{\ld,\mu}(q,t)$ is a polynomial in 
$q$ and $t$ with integral coefficients. 
} 
\noindent
(See Theorem 3.2).  
Let us recall (\cite{Ma1}, (VI.8.11)) that the {\it double Kostka coefficients} 
(or $(q,t)$-{\it Kostka coefficients}) 
$K_{\ld,\mu}(q,t)$ are defined via decomposition 
$$
J_{\mu}(x;q,t)=\sum_{\ld} K_{\ld,\mu}(q,t) S_{\ld}(x;t),
\tag{5}
$$
where $S_{\ld}(x;t)$ are the so-called {\it big Schur functions}
 (\cite{Ma1}, (III.4.5)). 
Theorem 2 gives a partial affirmative answer 
to the conjecture of Macdonald proposed in 
\cite{Ma1}, (VI.8.18?). 
\par
It also turns out that the raising operators $B_m$ and the lowering operators 
$A_m$ preserve the ring of symmetric functions $\K[x]^W$. 
On $\K[x]^W$, 
the action of these operators is described by the following $W$-invariant 
$q$-difference operators 
$$
\align
(B_m)_{\text{\rm sym}}
&
=\sum_{r=0}^m (-1)^r t^{\binom{r}{2}+(m-n+1)r}
\sum\Sb I\subset[1,n]\\ |I|=r\endSb
x_I e_{m-r}(x_{[1,n]\backslash I}) a_I(x) T_{q,x}^I, 
\tag{6}\\
(A_m)_{\text{\rm sym}}
&=\sum_{r=0}^m (-1)^r t^{\binom{r}{2}}
\sum\Sb I\subset[1,n]\\ |I|=r\endSb
x^{-1}_I e_{m-r}(x^{-1}_{[1,n]\backslash I}) a_I(x) T_{q,x}^I,
\endalign
$$
respectively (Theorem 7.1).  Here we used the abbreviated notations 
$$
x_I=\prod_{i\in I}x_i, \quad T_{q,x}^I=\prod_{i\in I}T_{q,x_i}, \quad
a_I(x)=\prod\Sb i\in I\\ j\not\in I\endSb \frac{tx_i-x_j}{x_i-x_j}, 
\tag{7}
$$ 
and, for a subset $J\subset[1,n]$, 
$e_{s}(x_J)$ denotes the elementary 
symmetric function of degree $s$ in the variables 
$(x_j)_{j\in J}$.
\par
After some preliminaries on Macdonald's $q$-difference operators and
the Dunkl operators (Sections 1 and 2), we formulate our main results in 
Section 3.  
Theorem 1 will be proved in Sections 5 and 6 by means of the 
{\it Mimachi basis} (defined in Section 4), 
which is a family of rational functions in two sets of variables that 
realizes some representation of the Hecke algebra.  
Explicit formulas (6) for the raising and the lowering 
$q$-difference operators will be determined in 
Section 7.  
\par
In the last section, we will give an application of our results to 
some combinatorial problem. 
We will introduce a double analogue of the multinomial coefficients 
in terms of the so-called modified Macdonald polynomials. 
The {\it modified Macdonald polynomials} $\widetilde{J}_\ld(x;q,t)$ 
are defined using the $\lambda$-ring notations 
as 
$$
\widetilde{J}_{\ld}(x;q,t) = J_\ld(\frac{x}{1-t};q,t). 
\tag{8}
$$
It is well-known (see e.g. \cite{GH}) that the double Kostka coefficients 
are characterized also as the transition coefficients between the modified 
Macdonald polynomials and the Schur functions: 
$$
\widetilde{J}_{\mu}(x;q,t)=\sum_{\mu} K_{\ld,\mu}(q,t) \, s_{\ld}(x). 
\tag{9}
$$
Let us introduce a family of polynomials $B_{\ld,\mu}(q,t)$ via decomposition
$$
\widetilde{J}_{\ld}(x;q,t)=\sum_{\mu} B_{\ld,\mu}(q,t) \, m_{\mu}(x)
\tag{10}
$$
in terms of the monomial symmetric functions.
The polynomiality of these coefficients follows from our Theorem 2. 
Note also that $B_{\ld,\mu}(q,t)=0$ unless $|\ld|=|\mu|$.
\Theorem{3}{
For any partitions $\ld$ and $\mu$ with $|\ld|=|\mu|$, we have 
\roster
\item \quad $B_{\ld,\mu}(q,t)\in \Z[q,t]$,
\item \quad $ B_{\ld,\mu}(1,1)=\left(\matrix
|\mu| \\ \mu_1,\mu_2,\cdots \endmatrix\right)$, 
\item \quad $ B_{(\ell),\mu}(q,t)=q^{n(\mu')}\left[\matrix
|\mu| \\ \mu_1,\mu_2,\cdots \endmatrix\right]_q$
\ \  if \ \  $\ell=|\mu|$,
\item \quad $ B_{\ld',\mu}(q,t)=q^{n(\ld')}t^{n(\ld)} B_{\ld,\mu}(t^{-1},q^{-1})$. 
\endroster
}
\noindent
It follows from Macdonald's conjecture \cite{Ma1}, (VI.8.18?) that
\proclaim{Conjecture 4?}  $B_{\ld,\mu}(q,t)\in \N[q,t]$ for any partitions 
$\ld$ and $\mu$. 
\endproclaim
\noindent
Hence, one can consider the polynomials $B_{\ld,\mu}(q,t)$ as a natural 
two-parameter deformation of the classical multinomial coefficients. 
\par\medpagebreak\par
The authors would like to express their thanks to Professor Katsuhisa Mimachi 
for valuable discussions. 
\par\medpagebreak
{\it Notes\,}: 
After we completed this paper, we found a direct proof for the 
fact that the $q$-difference operators (6) are raising operators 
for Macdonald polynomials. 
This method, without Dunkl operators, also provides an elementary 
proof of the integrality of double Kostka coefficients. 
For this direct approach, see our forthcoming paper \cite{KN}.  
\par\newpage
%
\section{\qquad Contents}
\par \qquad\S1:  Macdonald's $q$-difference operators.
\par \qquad\S2:  Affine Hecke algebras and the Dunkl operators.
\par \qquad\S3:  Raising operators and transition coefficients.
\par \qquad\S4:  Mimachi basis and a representation of the Hecke algebra.
\par \qquad\S5:  Action of $D_y(u)$ on $\Pi(x,y)$.
\par \qquad\S6:  Computation of the Dunkl operators acting on $\Pi(x,y)$.
\par \qquad\S7:  $q$-Difference raising operators. 
\par \qquad\S8:  A double analogue of the multinomial coefficients. 
\par\smallpagebreak
\par\qquad References
\par\newpage
\section{\S1:  Macdonald's $q$-difference operators}
%
In this section, we will make a brief review of some basic properties 
of the Macdonald polynomials (associated with the root system of type 
$A_{n-1}$, or the symmetric functions with two parameters) and 
the commuting family of $q$-difference operators
which have Macdonald polynomials as joint eigenfunctions.  
For details, see Macdonald's book \cite{Ma1}. 
\par\smallpagebreak\par
Let $\K=\Q(q,t)$ be the field of rational functions in two indeterminates 
$q$, $t$ 
and consider the ring $\K[x]=\K[x_1,\cdots,x_n]$ of polynomials in $n$ 
variables $x=(x_1,\cdots,x_n)$ with coefficients in $\K$.  
Under the natural action of the symmetric group $W=\frak{S}_n$ of 
degree $n$, the subring of all symmetric polynomials will be denoted 
by $\K[x]^{W}$. 
\par 
The {\it Macdonald polynomials} $P_{\ld}(x)=P_{\ld}(x;q,t)$ 
(associated with the root system of type $A_{n-1}$) 
are symmetric polynomials parametrized by the 
{\it partitions } $\ld=(\ld_1,\cdots,\ld_n)$
($\ld_i\in\Z$, $\ld_1\ge\cdots\ge\ld_n\ge0$).  
They form a $\K$-basis of the invariant ring $\K[x]^W$
and are characterized as the joint eigenfunctions 
of a commuting family of $q$-difference operators 
$\{D_x^{(r)}\}_{r=0}^n$. 
For each $r=0,1,\cdots,n$, the $q$-difference operator $D_x^{(r)}$ 
is defined by 
$$
D_x^{(r)}=\sum\Sb I\subset[1,n] \\ |I|=r \endSb t^{\binom{r}{2}}
\prod\Sb i\in I\\ j\not\in I\endSb 
\frac{tx_i-x_j}{x_i-x_j} \, \prod_{i\in I}T_{q,x_i}, 
\tag{1.1}
$$
where 
$T_{q,x_i}$ stands for the $q$-shift operator in the variable 
$x_i$ : 
$(T_{q,x_i}f)(x_1,\cdots,x_n)=f(x_1,\cdots,q x_i,\cdots,x_n)$. 
The summation in (1.1) is taken over all subsets $I$ of the interval 
$[1,n]=\{1,2,\cdots,n\}$ consisting of $r$ elements. 
Note that $D_x^{(0)}=1$ and 
$D_x^{(n)}=t^{\binom{n}{2}}T_{q,x_1}\cdots T_{q,x_n}$. 
Introducing a parameter $u$, we will use the generating function
$$
D_x(u)=\sum_{r=0}^n (-u)^r D_x^{(r)}
\tag{1.2} 
$$
of these operators $\{D_x^{(r)}\}_{r=0}^n$. 
Note that the operator $D_x(u)$ has the determinantal expression
$$
D_x(u)=\frac{1}{\Dt(x)}\det(x_j^{n-i}(1-ut^{n-i}T_{q,x_j}); 1\le i,j\le n),
\tag{1.3}
$$
where $\Dt(x)=\prod_{1\le i<j\le n}(x_i-x_j)$ is the difference product 
of $x_1,\cdots,x_n$. 
It is well known that the $q$-difference operators $D_x^{(r)}$ ($0\le r\le n$) 
commute with each other, or equivalently, $[D_x(u), D_x(v)]=0$. 
Furthermore the Macdonald polynomial $P_{\ld}(x)$ satisfies the 
$q$-difference equation
$$
D_x(u) P_{\ld}(x)=c^n_\ld(u) P_{\ld}(x), \quad\text{with}\quad 
c^n_\ld(u)=\prod_{i=1}^n (1-ut^{n-i}q^{\ld_i}), 
\tag{1.4}
$$
for each partition $\ld=(\ld_1,\cdots,\ld_n)$. 
Recall that each $P_\ld(x)$ can be written in the 
form
$$
P_\ld(x)=m_{\ld}(x) + \sum_{\mu<\ld} u_{\ld\mu} m_{\mu}(x)\quad(u_{\ld\mu}\in \K),
\tag{1.5}
$$
where, for each partition $\mu$, $m_{\mu}(x)$ stands for the monomial 
symmetric function of monomial type $\mu$, and $\le$ is the dominance 
order of partitions. 
The Macdonald polynomials $P_\ld(x)$ are determined uniquely by the conditions 
(1.4) and (1.5). 
\par
We recall here on the ``reproducing kernel'' of the Macdonald 
polynomials. 
Consider another set of variables $y=(y_1,\cdots,y_m)$ and assume that 
$m\le n$.  
We define the function $\Pi(x,y)=\Pi(x,y;q,t)$ by
$$
\Pi(x,y)=\prod\Sb 1\le i\le n\\ 1\le j\le m\endSb
\frac{(tx_iy_j;q)_\infty}{(x_iy_j;q)_\infty}, 
\tag{1.6}
$$
where $(x;q)_\infty=\prod_{k=0}^\infty (1-xq^k)$. 
The convergence of the infinite product above may be 
understood in the sense of formal power series
(or of absolute convergence assuming that 
$q$ is a complex variable with $|q|<1$). 
It is known that the function $\Pi(x,y)$ has the expression 
$$
\Pi(x,y)=\sum_{\ell(\ld)\le m} b_\ld P_\ld(x) P_\ld(y) \quad(b_\ld\in \K),
\tag{1.7}
$$
where the summation is taken over all partitions $\ld$ with length $\le m$, 
and each partition $\ld=(\ld_1,\cdots,\ld_m,0,\cdots,0)$ with 
$\ell(\ld)\le m$
is identified with the truncation $(\ld_1,\cdots,\ld_m)$ when it is used as the 
suffix for $P_\ld(y)$. 
The coefficients $b_\ld=b_\ld(q,t)$ in (1.7) are determined as 
$$
b_\ld=\prod_{s\in\ld} 
\frac{1-t^{\ell(s)+1}q^{a(s)}}{1-t^{\ell(s)}q^{a(s)+1}},
\tag{1.8}
$$
in terms of the leg-length $\ell(s)=\ld'_j-i$ and the arm-length 
$a(s)=\ld_i-j$ 
for a box $s=(i,j)$ in the Young diagram representing the 
partition $\ld$. 
\par
We remark that, by (1.4), expression (1.7) is equivalent to the formula
$$
D_x(u)\Pi(x,y)=(u;t)_{n-m} D_y(u t^{n-m}) \Pi(x,y). 
\tag{1.9}
$$
Since 
$$
D_x(u)=\sum_{I\subset[1,n]} (-u)^{|I|} t^{\binom{|I|}{2}}
\prod\Sb i\in I\\ j\not\in I\endSb
\frac{tx_i-x_j}{x_i-x_j} \, \prod_{i\in I}T_{q,x_i}, 
\tag{1.10}
$$
We have 
$$
D_x(u) \Pi(x,y) = \Pi(x,y) F(u;x,y),
\tag{1.11}
$$
where
$$
F(u;x,y)=\sum_{I\subset[1,n]} (-u)^{|I|} t^{\binom{|I|}{2}}
\prod\Sb i\in I\\ j\not\in I\endSb 
\frac{tx_i-x_j}{x_i-x_j} \, 
\prod\Sb i\in I\\ 1\le k\le m\endSb
\frac{1-x_iy_k}{1-tx_iy_k}.
\tag{1.12}
$$
Hence formula (1.9) is also equivalent to 
$$
F(u;x,y)=(u;t)_{n-m} F(ut^{n-m};y,x). 
\tag{1.13}
$$
(See also \cite{MN1}.)
\par\newpage
\section{\S2: Affine Hecke algebras and the Dunkl operators}
By the work of I.\,Cherednik \cite{C}, 
it is known that Macdonald's $q$-difference operators $\{D_x^{(r)}\}$ 
are reconstructed from the structure of affine Hecke algebras. 
We recall here how the Dunkl operators for Macdonald polynomials are 
defined, and how Macdonald's commuting family of $q$-difference operators 
are recovered from the Dunkl operators.  (See also I.G.\,Macdonald \cite{Ma2}, 
and A.A.\,Kirillov Jr. \cite{KJr}.) 
Our convention of the affine Hecke algebra and the definition of 
Dunkl operators are slightly different from the ones in the references 
cited above. 
\par\smallpagebreak\par
We denote by $P=\Z\ep_1\oplus\cdots\oplus\Z\ep_n$ the free $\Z$-module 
of rank $n$ with basis $\{\ep_i\}_{i=1}^n$
and take the canonical symmetric bilinear form 
$\scpr{\ }{}:P\times P \to \Z$ such that $\scpr{\ep_i}{\ep_j}=\dt_{ij}$
for each $1\le i,j\le n$.  
The elements of $P$ will be called (integral) {\it weights}, 
and a weight $\ld=\sum_{i=1}^n\ld_i\ep_i \in P$ will be 
identified freely with the $n$-tuple of integers 
$(\ld_1,\cdots,\ld_n)$.  
The action of the Weyl group $W=\frak{S}_n$ on $P$ will be fixed 
so that $w(\ep_i)=\ep_{w(i)}$ for $i=1,\cdots,n$,  
namely,  $w(\ld)_i=\ld_{w^{-1}(i)}$ for each $\ld\in P$ and 
$i=1,\cdots,n$. 
We will take the 
{\it simple roots} $\alpha_i=\ep_i-\ep_{i+1}$ and the 
{\it simple transpositions} $s_i=(i,i+1)$  for $i=1,\cdots,n-1$,
as usual. 
\par 
{}From this section on, we set $\tau_i=\tau_{x_i}=T_{q,x_i}$ ($i=1,\cdots,n$) 
to avoid the conflict with the notation of Hecke algebras. 
We use the notation of multi-indices both for the multiplication 
operators and for the $q$-shift operators: 
$$ 
x^\ld=x_1^{\ld_1}\cdots x_n^{\ld_n}\quad \text{and}\quad
\tau^\mu=\tau_1^{\mu_1}\cdots\tau_n^{\mu_n}
\tag{2.1}
$$ 
for each $\ld=(\ld_1,\cdots,\ld_n),\mu=(\mu_1,\cdots,\mu_n)\in P$. 
We denote by ${\Cal D}_{q,x}=\K(x)[\tau^{\pm 1}]$ the 
$\K$-subalgebra of $\operatorname{End}_\K(\K(x))$ generated by 
the multiplication by elements of $\K(x)$ and
the $q$-shift operators $\tau^\mu$ ($\mu\in P$): 
$$
{\Cal D}_{q,x}=\K(x)[\tau^{\pm 1}]=\bigoplus_{\mu\in P}\K(x) \tau^\mu.
\tag{2.2}
$$ 
Note that the $q$-shift operators $\tau^\mu$ act on $\K(x)$ 
as $\K$-algebra automorphisms and that 
the commutation relations between multiplication operators 
and $q$-shift operators are determined accordingly. 
In particular we have
$$
\tau^\mu x^\ld = q^{\scpr{\mu}{\ld}} x^\ld \tau^\mu\quad(\ld,\mu\in P).
\tag{2.3}
$$
Each element $w$ of the Weyl group $W=\frak{S}_n$ acts on $\K(x)$ 
as the $\K$-algebra automorphism of $\K(x)$ such that 
$w(x_i)=x_{w(i)}$ for $i=1,\cdots,n$.  
We denote by ${\Cal D}_{q,x}[W]$ the $\K$-algebra of 
$q$-difference operators involving permutations:  
$$
{\Cal D}_{q,x}[W]=\bigoplus_{w\in W} {\Cal D}_{q,x} w
=\bigoplus_{\mu\in P, w\in W}  \K(x) \tau^{\mu}w. 
\tag{2.4}
$$
Note that we have the commutation relations 
$$
w \,x^\ld = x^{w(\ld)} w \quad\text{and}\quad 
w \,\tau^\mu = \tau^{w(\mu)} w 
\tag{2.5}
$$
for all $\ld,\mu\in P$ and $w\in W$. 
\par
The $\K$-subalgebra $\K[\tau^{\pm1};W]$ of 
${\Cal D}_{q,x}[W]$, generated by the $q$-shift operators 
and the Weyl group, is isomorphic to the group ring 
$\K[\widetilde{W}]$ of the extended affine Weyl group 
$\widetilde{W}=P\rtimes W$. 
Let us describe the commutation relations of this algebra
in terms of generators. 
Firstly, the Weyl group $W$ is generated by the simple 
transpositions $s_i$ ($i=1,\cdots,n-1$).  
We define the element $s_0$ (corresponding to the affine 
simple root $\dt-\ep_1+\ep_n$) by
$$
s_0= s_{\ep_1-\ep_n} \tau_1\tau_n^{-1}.
\tag{2.6}
$$
These elements $s_0,s_1,\cdots,s_{n-1}$ generate the 
affine Weyl group $W^{\text{aff}}=Q^\vee\rtimes W$, $Q^\vee$ being 
the coroot lattice.  
Note that $W^{\text{aff}}=Q^\vee\rtimes W$ is a subgroup of 
our $\widetilde{W}=P\rtimes W$. 
The fundamental relations among $s_0,s_1,\cdots,s_{n-1}$ 
are given by 
$$
\alignat{2}
\text{(i)}&\quad s_i^2=1 &\quad &(i=0,1,\cdots,n-1),\\
\text{(ii)}&\quad s_is_j=s_js_i &\quad &(|i-j|\ge 2), \tag{2.7}\\
\text{(iii)}&\quad s_is_js_i=s_js_is_j &\quad & (|i-j|=1),
\endalignat
$$
if $n\ge 3$.  
Here the suffices for $s_0, s_1,\cdots,s_{n-1}$ are understood 
as elements of $\Z/n\Z$, and $|a|$ stands for the representative 
$r$ of the class $a+n\Z$ such that $0\le r<n$. 
If $n=2$, the fundamental relations are simply given by (2.7.i). 
In order to obtain the whole extended affine Weyl group 
$\widetilde{W}=P\rtimes W$, we need to adjoin an element,
denoted by $\omega$ below,  corresponding to the rotation of the 
Coxeter diagram. 
We set
$$
\omega=s_{n-1}s_{n-2}\cdots s_1 \tau_1=\tau_n s_{n-1}s_{n-2}\cdots s_1. 
\tag{2.8}
$$
As to this element $\omega$, we have the commutation relations 
$$
\text{(iv)}\quad\omega \, s_i= s_{i-1} \, \omega \quad(i=0,1,\cdots, n-1).
\tag{2.9}
$$
We remark that $\omega$ has the infinite order in our $\widetilde{W}$, and 
that $\omega^n$ coincides with the Euler operator $\tau_1\cdots\tau_n$. 
Summarizing, the extended affine Weyl group $\widetilde{W}=P\rtimes W$ 
is generated by $s_0,s_1,\cdots,s_{n-1}$ and $\omega$, and their 
fundamental relations are given by (i) -- (iv) in (2.7) and (2.9). 
Note also that the $q$-shift operators $\tau_1,\cdots,\tau_n$ are 
recovered by the formula
$$
\tau_i=s_i s_{i+1}\cdots s_{n-1} \,\omega\, s_1 \cdots s_{i-1} 
\tag{2.10}
$$
for $i=1,\cdots,n-1$. 
\par
One important fact is that the Hecke algebra $H(\widetilde{W})$ 
of the extended affine Weyl 
group $\widetilde{W}=P\rtimes W$ can be realized in the algebra 
${\Cal D}_{q,x}[W]$ of $q$-difference operators with permutations. 
We define the elements $T_i$ ($i=0,1,\cdots n-1$) in 
${\Cal D}_{q,x}[W]$ by
$$
T_i=t+\frac{1-tx_i/x_{i+1}}{1-x_i/x_{i+1}}(s_i -1)
\tag{2.11}
$$
for $i=1,\cdots, n-1$ and  
$$
T_0=t+\frac{1-tqx_n/x_1}{1-qx_n/x_1}(s_0 -1). 
\tag{2.12}
$$
Then the following relations are verified in ${\Cal D}_{q,x}[W]$: 
$$
\alignat{2}
\text{(i)}&\quad (T_i-t)(T_i+1)=0 &\quad &(i=0,1,\cdots,n-1),\\
\text{(ii)}&\quad T_iT_j=T_jT_i &\quad &(|i-j|\ge 2), \tag{2.13}\\
\text{(iii)}&\quad T_iT_jT_i=T_jT_iT_j &\quad & (|i-j|=1),\\
\text{(iv)}&\quad\omega \, T_i= T_{i-1} \, \omega & \quad & (i=0,1,\cdots, n-1),
\endalignat
$$
with indices understood as elements of $\Z/n\Z$.  
We denote by  $H(\widetilde{W})$ the subalgebra of ${\Cal D}_{q,x}[W]$ 
generated by $T_0,T_1,\cdots,T_{n-1}$ and $\omega^{\pm1}$.  
One can also show that (2.13) gives a complete list of the 
fundamental relations among the generators $T_0,T_1,\cdots,T_{n-1}$ 
and $\omega^{\pm1}$.
This is a realization of the extended affine Hecke algebra 
(defined by the generators $T_0,\cdots,T_{n-1}, \omega^{\pm 1}$ 
and the relations (i) -- (iv) above) 
in the algebra of $q$-difference operators with permutations. 
Note also that the $\K$-subalgebra of ${\Cal D}_{q,x}[W]$ generated 
by $T_1,\cdots,T_{n-1}$ is isomorphic to the Hecke algebra 
$H(\frak{S}_n) $ of the symmetric group. 
In what follows, 
we will also use the operators $\overline{T}_i=t^{-1}T_i $ 
($i=0,1,\cdots,n-1$), which satisfy the quadratic relations 
$ (\overline{T}_i-1)(\overline{T}_i+t^{-1})=0$. 
\par
We now define the {\it Dunkl operators} $Y_1,\cdots,Y_n$ and 
the {\it dual Dunkl operators} $Y^\ast_1,\cdots,Y^\ast_n$ in the
extended affine Hecke algebra $H(\widetilde{W})$.  
For each $i=1,\cdots,n$, we set
$$
\align
Y_i&=\overline{T}_i \overline{T}_{i+1}\cdots \overline{T}_{n-1}\,\omega\, 
\overline{T}_1^{-1} \cdots\overline{T}_{i-1}^{-1} \tag{2.14}\\
&=t^{-n+2i-1}\, 
T_i T_{i+1}\cdots T_{n-1}\,\omega\, T_1^{-1}\cdots T_{i-1}^{-1}
\endalign
$$
and
$$
\align
Y^\ast_i&=\overline{T}_i^{-1} \overline{T}_{i+1}^{-1}\cdots 
\overline{T}_{n-1}^{-1}\,\omega\, 
\overline{T}_1\cdots\overline{T}_{i-1} \tag{2.15}\\
&=t^{n-2i+1}
T_i^{-1} T_{i+1}^{-1}\cdots T_{n-1}^{-1}\,\omega\, T_1\cdots T_{i-1}. 
\endalign
$$
Note that 
$Y_1=\overline{T}_1 \overline{T}_{2}\cdots \overline{T}_{n-1}\,\omega$ 
and 
$Y_n=\omega\, 
\overline{T}_1^{-1} \cdots\overline{T}_{n-1}^{-1}$.
Comparing these formulas with (2.10), one sees that 
both $Y_i$ and $Y^\ast_i$ reduce to $\tau_i$ when $t\to1$.  
By (2.13), one can show that the Dunkl operators 
$Y_1,\cdots,Y_n$ commute with each other, and that
$$
H(\widetilde{W}) =\bigoplus_{w\in W}\K[Y^{\pm1}] T_w
=\bigoplus_{\mu\in P,\,w\in W}\K\,Y^\mu T_w,
\tag{2.16}
$$
where $Y^\mu=Y_1^{\mu_1}\cdots Y_n^{\mu_n}$. 
For each $w\in W$, $T_w$ and $\overline{T}_w$ is 
defined by
$$
T_w=T_{i_1}\cdots T_{i_p},\quad 
\overline{T}_w=\overline{T}_{i_1}\cdots 
\overline{T}_{i_p}=t^{-\ell(w)}T_w
\tag{2.17}
$$
by taking any reduced decomposition $w=s_{i_1}\cdots s_{i_p}$
$(1\le i_1,\cdots, i_p\le n-1)$;
these elements do not depend on the choice of reduced decompositions. 
The Dunkl operators satisfy the following commutation relations 
with $\overline{T}_1,\cdots,\overline{T}_{n-1}$:
$$
\align
&\overline{T}_{i} Y_{i+1} \overline{T}_i = Y_{i}, 
\quad\overline{T_i} Y_j= Y_j\overline{T_i} \quad\quad(j\ne i,i+1),
\tag{2.18}\\
&
\overline{T}_{i} Y^\ast_{i} \overline{T}_i = Y^\ast_{i+1}, 
\quad \overline{T_i} Y^\ast_j= Y^\ast_j\overline{T_i} \quad\quad(j\ne i,i+1),
\endalign
$$
for $i=1,\cdots,n-1$.
\par
Let us define a $\K$-algebra homomorphism 
$\eta : \K[\tau^{\pm 1}] \to \K[Y^{\pm 1}]$ by the substitution 
$\tau_i\mapsto t^{n-i} Y_i$ for $i=1,\cdots,n$. 
Then it is known that $\eta$ induces the isomorphism
$$
\eta :  \K[\tau^{\pm 1}]^W @>\sim>> {\Cal Z}H(\widetilde{W}) 
\tag{2.19}
$$
from the invariant ring $\K[\tau^{\pm 1}]^W$ onto the center 
of $H(\widetilde{W})$ (Bernstein's theorem).
This theorem also implies that, 
for any $f=f(\tau)\in \K[\tau^{\pm 1}]^W$, 
the operator 
$\eta(f)=f(t^{n-1}Y_1,t^{n-2}Y_2,\cdots,Y_n)\in{\Cal D}_{q,x}[W]$
is $W$-invariant, hence preserves the $\K$-algebra 
$\K[x]^W$ of symmetric polynomials. 
In this way, the center ${\Cal Z}H(\widetilde{W})$ of the extended 
affine Hecke algebra provides a commuting family of 
$q$-difference operators acting on $\K[x]^W$.  
It turns out also that this family of operators is diagonalized
simultaneously by the Macdonald polynomials $P_\ld(x)$. 
In fact we have 
$$
f(t^{n-1}Y_1,t^{n-2}Y_2,\cdots,Y_n) P_\ld(x)
=f(t^{n-1}q^{\ld_1},t^{n-2}q^{\ld_2},\cdots,q^{\ld_n}) P_\ld(x),
\tag{2.20}
$$
for any partition $\ld=(\ld_1,\cdots,\ld_n)$.  
In particular, the $q$-difference operator  $D_x(u)$ of Macdonald 
can be recovered as the restriction of an operator in 
${\Cal Z}H(\widetilde{W})$.  
Namely we have 
$$
(1-ut^{n-1}Y_1)(1-ut^{n-2}Y_2)\cdots(1-uY_n) \vert_{\K[x]^W}
=D_x(u).
\tag{2.21}
$$
\par
We remark that the dual Dunkl operators $Y_1^\ast,\cdots,Y_n^\ast$ defined 
by (2.15) also have properties similar to $Y_1,\cdots,Y_n$ in relation to 
Macdonald polynomials.  
They commute with each other, and are related with Macdonald's 
$q$-difference operator through the formula
$$
(1-uY^\ast_1)(1-utY^\ast_2)\cdots(1-ut^{n-1}Y^\ast_n) \vert_{\K[x]^W}
=D_x(u).
\tag{2.22}
$$
(See Section 4.) 
%
\section{\S3: Raising operators and transition coefficients}
In this section we introduce certain raising and lowering operators 
for Macdonald polynomials.  
After stating our main result (Theorem 3.1), 
we discuss some of its consequences, 
including the polynomiality and the integrality of 
transition coefficients related to Macdonald polynomials.   
At the end of this section, we formulate the key lemma 
for the proof of Theorem 3.1, and give a direct proof for a special case, 
to show some of the ideas in our proof of the general case in 
Sections 4--6. 
\par\medpagebreak
By means of the Dunkl and the dual Dunkl operators, we can construct 
certain raising and lowering operators for Macdonald polynomials. 
Using the Dunkl operators $Y_1,\cdots, Y_n$ of (2.14), 
we define the operator $B^x_m \in {\Cal D}_{q,x}[W]$ by
$$
B^x_m=\sum_{1\le k_1<\cdots<k_m\le n}\,
x_{k_1}\cdots x_{k_m}
(1-t^m Y_{k_1})(1-t^{m-1}Y_{k_2})\cdots(1-tY_{k_m}).  
\tag{3.1}
$$
for each $m=1,2,\cdots,n$. 
Similarly we define the operator  $A^x_m\in {\Cal D}_{q,x}[W]$ by
$$
A^x_m=\sum_{1\le k_1<\cdots<k_m\le n}\,
\frac{1}{x_{k_1}\cdots x_{k_m}}
(1-Y_{k_1}^\ast)(1-t Y_{k_2}^\ast)\cdots(1-t^{m-1}Y_{k_m}^\ast), 
\tag{3.2}
$$
where the $Y^\ast_1, \cdots,Y^\ast_n$ are the dual Dunkl operators defined 
in (2.15). 
\Theorem{3.1}{
{\rm (1)} The operators $B^x_m\in{\Cal D}_{q,x}[W]$ $(m=1,\cdots,n)$ 
are raising operators 
for the Macdonald polynomials such that
$$
B^x_m \, P_\ld(x)=
\prod_{i=1}^m(1-t^{m-i+1}q^{\ld_i})
P_{\ld+(1^m)}(x)
\tag{3.3}
$$
for all partition $\ld$ with $\ell(\ld)\le m$, where 
$(1^m)=\ep_1+\cdots+\ep_m$. 
\newline
{\rm (2)} The operators $A^x_m\in{\Cal D}_{q,x}[W]$ $(m=1,\cdots,n)$ 
are lowering operators 
for the Macdonald polynomials such that
$$ 
A^x_m \, P_\ld(x)=
\prod_{i=1}^m
\frac{(1-t^{m-i}q^{\ld_i})(1-t^{n-i+1}q^{\ld_i-1})}{(1-t^{m-i+1}q^{\ld_i-1})}
P_{\ld-(1^m)}(x)
\tag{3.4}
$$
if $\ell(\ld)\le m$. 
In particular, one has $A^x_m \, P_\ld(x)=0$ if $\ell(\ld)<m$. 
}
\noindent
The proof of Theorem 3.1 will be given later in Sections 4, 5 and 6.
On the $\K$-algebra $\K[x]^W$ of symmetric polynomials, the operators 
$B^x_m$ and $A^x_m$ act as $q$-difference operators. 
Explicit formulas for these $q$-difference operators will be given 
in Section 7. 
In this section, we will discuss some of the consequences of 
Theorem 3.1. 
\par 
{}From Theorem 3.1, it follows that the Macdonald polynomial $P_\ld(x)$ 
for any partition $\ld$ can be obtained from the constant function $P_0(x)=1$ 
by an iterated application of the raising operators $B^x_m$. 
Namely we have
$$
P_\ld(x)=\text{const.\,} 
(B^x_n)^{\ld_n}(B^x_{n-1})^{\ld_{n-1}-\ld_n}\cdots (B^x_1)^{\ld_1-\ld_2}(1). 
\tag{3.5}
$$
In other words, 
$$
P_\ld(x)=\text{const.\,} 
B^x_{\mu_1}B^x_{\mu_2}\cdots B^x_{\mu_s}(1), 
\tag{3.6}
$$
where $\mu=(\mu_1,\cdots,\mu_s)$ 
($\mu_1\ge\cdots\ge\mu_s>0$)
is the conjugate partition $\ld'$ of $\ld$. 
If we take the normalization 
$$
J_\ld(x)=c_\ld P_\ld(x) ,\quad c_\ld=\prod_{s\in\ld} (1-t^{\ell(s)+1}q^{a(s)}), 
\tag{3.7}
$$
of the Macdonald polynomials (\cite{Ma1}), the action of the operator 
$B^x_m$ and $A^x_m$ are given as follows:
$$
\align
&B^x_m J_\ld(x)=J_{\ld+(1^m)}(x), \tag{3.8}\\
&A^x_m J_\ld(x)=\prod_{i=1}^m(1-t^{m-i}q^{\ld_i})(1-t^{n-i+1}q^{\ld_i-1})
\, J_{\ld-(1^m)}(x)
\endalign
$$
for $\ell(\ld)\le m$.  
In particular we have
$$
J_\ld(x)
=(B^x_n)^{\ld_n}(B^x_{n-1})^{\ld_{n-1}-\ld_n}\cdots (B^x_1)^{\ld_1-\ld_2}(1). 
\tag{3.9}
$$
Namely, the ``const.'' in (3.5) and (3.6) is precisely 
the reciprocal of the $c_\ld$ mentioned above.  
\par\medpagebreak
We can make use of expression (3.9) to study transition coefficients 
for Macdonald polynomials $J_\ld(x)=J_\ld(x;q,t)$. 
By (1.5) and (3.7), the ``integral form" $J_\ld(x;q,t)$ 
has an expression
$$
J_\ld(x;q,t)=c_\ld(q,t) \sum_{\mu\le\ld}  u_{\ld\mu}(q,t)\, m_\mu(x),
\tag{3.10}
$$
where $c_{\ld}(q,t)=c_\ld\in\Z[q,t]$ as in (3.7)
and $u_{\ld\mu}(q,t)\in\K=\Q(q,t)$ 
($u_{\ld\ld}(q,t)=1$). 
Following Macdonald \cite{Ma1}, Chapter VI, 
let us define the {\it double Kostka coefficients} $K_{\ld,\mu}(q,t)$ 
via decomposition
$$
J_{\mu}(x;q,t)=\sum_{\ld} K_{\ld,\mu}(q,t) S_{\ld}(x;t),
\tag{3.11}
$$
in terms of the {\it big Schur functions} $S_\ld(x;t)$
(by taking the stable limit as $n\to \infty$, to be more precise). 
For the definition of $S_{\ld}(x;t)$, we refer to \cite{Ma1}, (III.4.5). 
By means of expression (3.9), we can prove that 
$c_\ld(q,t)u_{\ld\mu}(q,t)$ and $K_{\ld\mu}(q,t)$
are in the polynomial ring $\Z[q,t]$ with integral coefficients,
for any partitions $\ld$ and $\mu$. 
\Theorem{3.2}{
{\rm (1)} 
For each partition $\ld$, the Macdonald polynomial $J_\ld(x;q,t)$ is  
expressed as a linear combination of monomial symmetric functions  
with coefficients in the ring $\Z[q,t]$. 
\newline
{\rm (2)} For any partitions $\ld$ and $\mu$, the double Kostka coefficient
$K_{\ld,\mu}(q,t)$ is a polynomial in $q$ and $t$ with integral coefficients. 
}
\Proof{
Since the divided difference operators 
$\frac{tx_i-x_{i+1}}{x_i-x_{i+1}}(s_i-1)$ 
($i=1,\cdots,n-1$) preserve the ring $\Z[t,t^{-1}][x]$,  
so do the operators $T_i$ as well as $T_i^{-1}$.  
This implies that the Dunkl operators 
$Y_1,\cdots,Y_n$, and accordingly the raising operators $B^x_m$
preserve the ring $\Z[q;t,t^{-1}][x]^W$.  
Hence we have $J_\ld(x;q,t)\in\Z[q;t,t^{-1}][x]^W$. 
Note that this statement is valid for any $n$. 
Hence, also in infinite variables, 
$J_\ld(x_1,x_2,\cdots;q,t)$ is a linear combination of 
monomial symmetric functions with coefficients in $\Z[q;t,t^{-1}]$.
We can now apply the duality of $J_\ld(x_1,x_2,\cdots;q,t)$ and 
$J_{\ld'}(x_1,x_2,\cdots;t,q)$ with respect to $q$ and $t$ in 
\cite{Ma1}, (VI.8.6), 
and conclude that $J_\ld(x_1,x_2,\cdots;q,t)$ is  
a linear combination of monomial symmetric functions with coefficients 
regular at $t=0$, since the involution $\omega_{q,t}$ does not give rise to 
any singularity at $t=0$ or $q=0$.  
Hence the coefficients must belong eventually to $\Z[q,t]$.  
This also implies $J_\ld(x;q,t)\in \Z[q,t][x]^W$, 
which proves Statement (1). 
\newline 
As in \cite{Ma1}, p.241, the transition coefficients 
between monomial symmetric functions and big Schur functions 
have the form $p(t)/q(t)$, where $p(t), q(t)\in\Z[t]$ and $q(0)=1$. 
Note in particular that they belong to the ring $\Q[t]_{(t)}$ of 
rational functions in $t$, regular at $t=0$. 
Combining this fact with Statement (1) of Theorem, we see that each 
double Kostka coefficients can be written as a finite sum of the form
$$
K_{\ld,\mu}(q,t)=\sum_{k\ge 0} \frac{p^{(k)}_{\ld\mu}(t)}{q^{(k)}_{\ld\mu}(t)} q^k,
\tag{3.12}
$$
where all $p^{(k)}_{\ld\mu}(t)$ and $q^{(k)}_{\ld\mu}(t)$ belong to $\Z[t]$ and 
$q^{(k)}_{\ld\mu}(0)=1$. 
In particular we have $K_{\ld,\mu}(q,t)\in \Q[t]_{(t)}[q]$.  
We can now apply the duality with respect to $q$ and $t$ again, 
between $J_{\ld}(x;q,t)$ and $J_{\ld'}(x;t,q)$
and between $S_{\ld}(x;t)$ and $S_{\ld'}(x;q)$, to conclude 
$K_{\ld,\mu}(q,t)=K_{\ld',\mu'}(t,q)$ (\cite{Ma1},(VI.8.15)).  
Hence, we have 
$K_{\ld,\mu}(q,t)\in \Q[t]_{(t)}[q]\cap  \Q[q]_{(q)}[t]$. 
By using Taylor expansions at $t=q=0$, one can easily see 
that the intersection of the two subalgebras $\Q[t]_{(t)}[q]$ and 
$\Q[q]_{(q)}[t]$ coincides precisely with $\Q[q,t]$.  
Hence we have $K_{\ld,\mu}(q,t)\in \Q[q,t]$. 
In expression (3.12), it means that, for any $k$, 
$p^{(k)}_{\ld\mu}(t)/q^{(k)}_{\ld\mu}(t)\,\in\Q[t]$, 
i.e., $q^{(k)}_{\ld\mu}(t)$ divides $p^{(k)}_{\ld\mu}(t)$. 
Since $q^{(k)}_{\ld\mu}(0)=1$, it follows that 
$p^{(k)}_{\ld\mu}(t)/q^{(k)}_{\ld\mu}(t)\,\in\Z[t]$ 
for all $k$. 
Namely we have $K_{\ld,\mu}(q,t)\in \Z[q,t]$. 
}
\par\medpagebreak\par
Our raising operators $B^x_m$ ($m=0,1,\cdots,n$) can be regarded as 
a natural $q$-analogue of those introduced by Lapointe-Vinet \cite{LV1, LV2}
for Jack polynomials, although the limits of our $B^x_m$ does not 
give exactly their operators. 
For the comparison with the operators of Lapointe-Vinet, we look 
at the ``quasi-classical limits'' of our operators as $q\to 1$. 
Introducing the parameter $\beta$, let us take the limit as $q\to 1$ with 
rescaling $t=q^\beta$.  
Then our Dunkl operators $Y_k$ ($k=1,\cdots,n$) 
have the following quasi-classical limits:
$$
{\Cal D}_k=\lim_{q\to 1}\frac{1-Y_k}{1-q}
=x_k\frac{\partial\ \ }{\partial x_k} 
+\beta\sum_{\alpha>0}
\frac{\scpr{\alpha}{\ep_k}}{1-x^\alpha}(s_{\alpha}-1),
\tag{3.13}
$$
summed over all positive roots $\alpha=\ep_i-\ep_j$ ($i<j$). 
Note that these Dunkl operators commute with each other, 
i.e., $[{\Cal D}_i, {\Cal D}_j]=0$.
With these ${\Cal D_k}$, the quasi-classical limit 
$\lim_{q\to 1}(B^x_m/(1-q)^m)$ 
of the raising operator
$B^x_m$ is given by
$$
\sum_{k_1<\cdots<k_m}
x_{k_1}x_{k_2}\cdots x_{k_m}
({\Cal D}_{k_1}+m\beta)({\Cal D}_{k_2}+(m-1)\beta)
\cdots({\Cal D}_{k_m}+\beta). 
\tag{3.14}
$$
Similarly, the quasi-classical limits of the lowering operators 
$A^x_m$ are given by
$$
\sum_{k_1<\cdots<k_m}
\frac{1}{x_{k_1}x_{k_2}\cdots x_{k_m}}
({\Cal D}^\ast_{k_1})({\Cal D}^\ast_{k_2}+\beta)
\cdots({\Cal D}^\ast_{k_m}+(m-1)\beta),
\tag{3.15}
$$
where 
$$
{\Cal D}^\ast_k=\lim_{q\to 1}\frac{1-Y^\ast_k}{1-q}
=x_k\frac{\partial\ \ }{\partial x_k} 
-\beta\sum_{\alpha>0}
\frac{\scpr{\alpha}{\ep_k}}{1-x^{-\alpha}}(s_{\alpha}-1),
\tag{3.16}
$$
%
\par\medpagebreak\par
Theorem 3.1 will be proved by investigating the action of the operators 
$B^x_m$ and $A^x_m$ ($m\le n$)
on the function $\Pi(x,y)$ with the auxiliary variables 
$y=(y_1,\cdots,y_m)$.  
In fact we will prove the following key lemma in Sections 5 and 6.
\Lemma{3.3}{
The operators $B_m^x$ and $A_m^x$ act on $\Pi(x,y)$ for $y=(y_1,\cdots,y_m)$ 
as follows:
$$
\align
&B_m^x \, \Pi(x,y)=\frac{1}{y_1\cdots y_m}D_y(1)\,\Pi(x,y),\tag{3.17}\\
&A_m^x \, \Pi(x,y)=y_1\cdots y_m D_y(t^{n-m+1})\,\Pi(x,y). \tag{3.18}
\endalign
$$
}
\noindent
It is easy to see that Lemma 3.3 implies Theorem 3.1.  
In fact, we have 
$$
\align
y_1\cdots y_m \, B^x_m \, \Pi(x,y) 
&=\sum_{\ell(\ld)\le m} b_\ld B^x_m (P_\ld(x)) \, y_1\cdots y_m \, P_\ld(y)
\tag{3.19}\\
&=\sum_{\ell(\ld)\le m} b_\ld B^x_m (P_\ld(x)) \,P_{\ld+(1^m)}(y).
\endalign
$$
As to the action of $D_y(1)$, we have 
$$
\align
&D_y(1)  \, \Pi(x,y) =\sum_{\ell(\ld)\le m} b_\ld \, P_\ld(x) \, 
c^m_{\ld}(1)\,P_\ld(y) \tag{3.20}\\
&\quad=\sum_{\ell(\ld)\le m} 
c^m_{\ld+(1^m)}(1) \, b_{\ld+(1^m)} \, P_{\ld+(1^m)}(x) \, P_{\ld+(1^m)}(y),
\endalign
$$
since 
$c^m_{\ld}(1)=(1-t^{m-1}q^{\ld_1})\cdots(1-tq^{\ld_{m-1}})(1-q^{\ld_m})=0$ 
if $\ld_m=0$. 
Since (3.19) and (3.20) are equal to each other by (3.17), 
we obtain 
$$
b_\ld\, B^x_m (P_\ld(x)) = c^m_{\ld+(1^m)}(1) \, b_{\ld+(1^m)}
\, P_{\ld+(1^m)}(x), 
\tag{3.21}
$$
by comparing the coefficients of $P_{\ld+(1^m)}(y)$. 
This proves statement (1) of Theorem 3.1. 
A similar computation based on (3.18) shows that
$$
b_\ld\, A^x_m (P_\ld(x)) = 
c^m_{\ld-(1^m)}(t^{n-m+1}) \, b_{\ld-(1^m)}
\, P_{\ld-(1^m)}(x), 
\tag{3.22}
$$
which proves statement (2) of Theorem 3.1. 
\Remark{3.4}{
It follows from Theorem 3.1 that, if $\ld$ is a partition such that
$\ell(\ld)\le m\le n$, then 
$$
[D_x(u), B_m] J_\ld(x) = (c^n_{\ld+(1^m)}(u)-c^n_{\ld}(u)) B_m J_\ld(x). 
\tag{3.23}
$$
In the case of Jack polynomials, the last formula corresponds to 
Proposition 4.1 of Lapointe-Vinet \cite{LV2} and appears to be the main 
step of their proof of a Rodrigues type formula for Jack polynomials \cite{LV1}.
It would be interesting to find a direct proof of formula (3.23). 
}
\par\medpagebreak\par
Before the proof of Lemma 3.3, we will give a proof of formula (3.17) 
for the case $m=1$, to show some of the ideas of our proof. 
\par
We have to compute 
$$
B^x_1 \, \, \Pi(x,y)= \sum_{k=1}^n x_k (1-tY_k) 
\prod_{i=1}^n \frac{(tx_iy;q)_\infty}{(x_iy;q)_\infty}. 
\tag{3.24}
$$
Since $\Pi(x,y)$ is symmetric in $x=(x_1,\cdots,x_n)$,  
we have
$$
\align
Y_k \, \Pi(x,y) &= 
\overline{T}_k\,\overline{T}_{k+1}\cdots\overline{T}_{n-1}
\tau_{n} \Pi(x,y)  
\tag{3.25}\\
&
=\Pi(x,y) \, 
\overline{T}_k\,\overline{T}_{k+1}\cdots\overline{T}_{n-1}
\left(\frac{1-x_ny}{1-tx_ny}\right)\\
&=\Pi(x,y) \, 
\overline{T}_k\,\overline{T}_{k+1}\cdots\overline{T}_{n-1}
t^{-1}\left(1+\frac{t-1}{1-tx_ny}\right). 
\endalign
$$
Hence
$$
(1-tY_k) \, \Pi(x,y) 
=\Pi(x,y) \, 
\overline{T}_k\,\overline{T}_{k+1}\cdots\overline{T}_{n-1}
\left(\frac{1-t}{1-tx_ny}\right). 
\tag{3.26}
$$
The key identity in our computation is 
$$
\overline{T}_k \left(\frac{1}{1-tx_{k+1}y}\right)=
\frac{1}{1-tx_ky}\frac{1-x_{k+1}y}{1-tx_{k+1}y}\quad(k=1,\cdots,n-1). 
\tag{3.27}
$$
Using this identity repeatedly, we get
$$
\overline{T}_k\,\overline{T}_{k+1}\cdots\overline{T}_{n-1}
\left(\frac{1}{1-tx_ny}\right)
=\frac{1}{1-tx_ky} \prod_{k<i\le n}\frac{1-x_iy}{1-tx_iy}.
\tag{3.28}
$$
This implies 
$$
(1-tY_k) \, \Pi(x,y) =
\Pi(x,y) 
\left(\frac{1-t}{1-tx_ky}\prod_{k<i\le n}\frac{1-x_iy}{1-tx_iy}\right).
\tag{3.29}
$$
We now use the formula
$$
\frac{(1-t)x_ky}{1-tx_ky}=1-\frac{1-x_ky}{1-tx_ky}
\tag{3.30}
$$
to compute the summation 
$$
\align
\sum_{k=1}^n \frac{(1-t)x_ky}{1-tx_ky}\prod_{k<i\le n}\frac{1-x_iy}{1-tx_iy} 
&=\sum_{k=1}^n \left(\prod_{k<i\le n}\frac{1-x_iy}{1-tx_iy}
-\prod_{k\le i\le n}\frac{1-x_iy}{1-tx_iy}\right) \tag{3.31}\\
&=1-\prod_{1\le i\le n}\frac{1-x_iy}{1-tx_iy}.
\endalign
$$
Hence we obtain
$$
\align
y \sum_{k=1}^n x_k(1-tY_k) \, \Pi(x,y) 
&=\Pi(x,y) \left(1-\prod_{1\le i\le n}
\frac{1-x_iy}{1-tx_iy}\right)\tag{3.32}\\
&=(1-\tau_y)\Pi(x,y). 
\endalign
$$
This gives formula (3.17) of Lemma 3.3 for $m=1$. 
A similar computation can be done to prove (3.18) for $m=1$. 
\par\medpagebreak\par
For the proof of Lemma 3.3, we will make use of a certain class of 
rational functions in $(x,y)$  which is related to a systematic 
generalization of formula (3.28) above. 
Such rational functions, which we call {\it Mimachi basis} below, 
were proposed by \cite{Mi} in the study of integral representations 
of $q$-KZ equations. 
It actually gives a realization of some representations of the Hecke 
algebra $H(W)$ as is clarified in \cite{MN2}. 
\par
The proof of Lemma 3.3 will be divided into two parts.  
After reformulating the Mimachi basis in Section 4 
so as to fit for our purpose, 
we will describe in Section 5 the action of the operator $D_y(u)$ on 
$\Pi(x,y)$, in terms of the Mimachi basis. 
On the other hand, we will analyze in Section 6 the action of the operators 
$B^x_m$ and $A^x_m$ on $\Pi(x,y)$ by computing the action of Dunkl 
operators explicitly to establish the formulas (3.17) and (3.18).
%
\section{\S4:  Mimachi basis and a representation of the Hecke algebra}
In this section, we will formulate the Mimachi basis in our context. 
Keeping the notations of the previous section,
we work with the field $\K(x,y)$ of rational functions in 
$x=(x_1,\cdots,x_n)$ and $y=(y_1,\cdots,y_m)$. 
For the moment, we do {\it not} assume that $m\le n$.  
\par\medpagebreak\par
Note first that the field $\K(x,y)$ has a natural structure of left module over the 
Hecke algebra $H(W)=H(\frak{S}^x_n)$ defined through the operators
$$
T^x_i=t+\frac{1-tx_i/x_{i+1}}{1-x_i/x_{i+1}}(s_i-1)\quad (i=1,\cdots,n-1)
\tag{4.1}
$$
as in Section 2. 
We use the superscript $x$ to remind that these operators are realized 
as operators with respect to the variables $x=(x_1,\cdots,x_n)$.
As to this $H(\frak{S}^x_n)$-module structure of $\K(x,y)$, 
we construct some explicit finite dimensional subrepresentations 
of $\K(x,y)$.  
\par
Let us introducing some notations.  
Let $I\subset[1,n]$ and $K\subset[1,m]$ be 
subsets of indices for the variables $x$ and $y$, respectively, 
and assume that $|I|=|K|=r$ ($0\le r\le m\wedge n$). 
By writing these sets as $I=\{i_1<\cdots<i_r\}$
and $K=\{k_1<\cdots<k_r\}$,  we define a rational function 
$h^{I}_{K}=h^{I}_{K}(x,y)\in \K(x,y)$ as follows:
$$
\align
&h^{I}_{K}(x,y)\tag{4.2}\\
&=\sum_{\sigma\in\frak{S}_r}
\prod_{1\le\mu\le r}\left(\frac{t-1}{1-tx_{i_\mu}y_{k_{\sigma(\mu)}}}
\prod_{i_\mu<j\le n}
\frac{t(1-x_jy_{k_{\sigma(\mu)}})}{1-tx_jy_{k_{\sigma(\mu)}}}\right)
\prod_{1\le\mu<\nu\le r} 
\frac{ty_{k_{\sigma(\mu)}}-y_{k_{\sigma(\nu)}}}
{y_{k_{\sigma(\mu)}}-y_{k_{\sigma(\nu)}}}. 
\endalign
$$
Note that each $h^{I}_{K}(x,y)$ is a symmetric function in the variables 
$(y_{k_1},\cdots,y_{k_r})$, 
while, as a function of $x$, it strongly depends on the inclusion 
$I\subset[1,n]$.
\par
For a fixed $K\subset[1,m]$ with $|K|=r$, 
let $V_{n,r;K}$ be the $\K$-vector subspace of 
$\K(x,y)$ spanned by the rational functions $h^{I}_{K}(x,y)$ 
($I\subset[1,n], \,|I|=r $): 
$$
V_{n,r;K} = \sum\Sb I\subset[1,n]\\ |I|=r\endSb \K \, h^{I}_{K}(x,y)\subset \K(x,y). 
\tag{4.3}
$$
\Theorem{4.1}{
For each $K\subset[1,m]$ with $|K|=r$ $(0\le r\le m\wedge n)$,  
the vector subspace $V_{n,r;\,K}$ is an $H(\frak{S}^x_n)$-submodule 
of $\K(x,y)$.  
The rational functions $\{h^{I}_{K}(x,y);\,I\subset[1,n], \ |I|=r \}$
form a $\K$-basis of $V_{n,r,K}$; hence $\dim_\K V_{n,r; K}=\binom{n}{r}$.   
Furthermore, for each $i=1,\cdots,n-1$, 
the action of the operator $T^x_i$ on
$h^{I}_{K}=h^{I}_{K}(x,y)$  is described as follows:
$$
\alignat{2}
\text{\rm (i)}&\quad T^x_i \, h^{I}_{K} =t \,h^{I}_{K} 
\qquad\quad(i,i+1\in I\quad  &\text{or}\quad i,i+1\not\in I),\\ 
\text{\rm (ii)}&\quad T^x_i \, h^{I}_{K}= h^{s_i(I)}_{K} 
& (i\not\in I, i+1\in I),\tag{4.4}\\
\text{\rm (iii)}&\quad T^x_i \, h^{I}_{K}=t\,h^{s_i(I)}_{K}+(t-1)h^{I}_{K} 
&(i\in I, i+1\not\in I).
\endalignat
$$
}
\noindent
Note that formula (iii) of (4.4) is equivalent to 
$(T^x_i)^{-1} h^{I}_{K}=h^{s_i(I)}_{K}$. 
We will not prove this theorem here since it is a reformulation of 
a part of the results of \cite{MN2}. 
\par\medpagebreak\par
{}From formula (4.4), it turns out that, for each $K$, 
the $H(\frak{S}^x_n)$-module $V_{n,r;K}$ is isomorphic 
to the induced representation $\operatorname{Ind}^{H(\frak{S}_n)}
(\operatorname{triv}_{H(\frak{S}_{n-r}\times\frak{S}_{r})})$
of the trivial representation of $H(\frak{S}_{n-r}\times\frak{S}_{r})$. 
\par
For each subset $I\subset[1,n]$ with $|I|=r$, take the 
permutation $w_I\in\frak{S}_{r}$ defined by
$$
w_I=\pmatrix 
1&\cdots & n-r & n-r+1 & \cdots & n \\
j_1&\cdots &j_{n-r} & i_1 &\cdots & i_r
\endpmatrix,
\tag{4.5}
$$
where $I=\{i_1<\cdots<i_r\}$ and 
$[1,n]\backslash I=\{j_1<\cdots<j_{n-r}\}$.
Note that a reduced decomposition of $w_I$ is given by 
$$
w_I=s_{i_r}\cdots s_{n-1} s_{i_{r-1}}\cdots 
s_{n-2} \cdots s_{i_1}\cdots s_{n-r},
\tag{4.6}
$$
and that the length of $w_I$ is determined as
$$
\ell(w_I)=\sum_{\nu=1}^r ( n-r+\nu - i_\nu)= rn-\binom{r}{2}-\sum_{i\in I} i. 
\tag{4.7}
$$
{}From (4.4) it follows that 
$$
h^{I}_{K}(x,y)=T^x_{w_I} \, h^{[n-r+1,n]}_{K}(x,y), 
\tag{4.8}
$$
since $T^x_{w_I}$ has the expression
$$
T^x_{w_I}=T_{i_r}\cdots T_{n-1} T_{i_{r-1}}\cdots 
T_{n-2} \cdots T_{i_1}\cdots T_{n-r}. 
\tag{4.9}
$$ 
\par
Recall that each permutation $w\in\frak{S}_n$ can be uniquely 
written in the form
$w= w_I w' w''$ with some 
$I\subset[1,n]$ ($|I|=r$), $w'\in\frak{S}_{n-r}$, $w''\in\frak{S}_r$
and that $\ell(w)=\ell(w_I)+\ell(w')+\ell(w'')$.  
Hence the elements $T_{w_I}$ ($I\subset[1,n], |I|=r$) 
form a free basis of $H(\frak{S}_n)$ as a right 
$H(\frak{S}_{n-r}\times\frak{S}_r)$-module.
This shows that the $H(\frak{S}_n)$-module $V_{n,r;K}$ is isomorphic 
to the induced representation mentioned above. 
Note that the canonical generator of the induced representation 
corresponds to the rational function $h^{[n-r+1,n]}_{K}$. 
We will call the rational functions 
$\{h^{I}_{K}(x,y); I\subset[1,n],\ |I|=r\}$ 
the {\it Mimachi basis} of 
the representation $V_{n,r;K}=H(\frak{S}^x_n) \, h^{[n-r+1,n]}_{K}$. 
\par
These properties of the Mimachi basis will be used in the Section 6. 
\par\medpagebreak\par
For the proof of Lemma 3.3, we will need a {\it dual version} of the Mimachi 
basis as well. 
Let $\iota : \K(x,y)\to\K(x,y)$ be the involutive $\Q$-{\it algebra} automorphism 
determined by 
$$
\iota(q)=q^{-1}, \ \iota(t)=t^{-1}, \ \iota(x_i)=x_i^{-1}\ \ , \ \iota(y_k)=y_k^{-1},
\tag{4.10}
$$
for $i=1,\cdots,n$ and $k=1,\cdots,m$.  
Note that $\iota$ is {\it not} an anti-automorphism. 
It is easily checked that
$$
\iota(w f)=w \iota(f) \quad (w\in W=\frak{S}^x_n),\quad 
\iota(\tau_x^\mu f) =\tau_x^\mu \iota(f)\quad (\mu\in P),
\tag{4.11}
$$
for all $f\in\K(x,y)$. 
Hence the involution $\iota$ intertwines the action of the extended affine Weyl group
$\widetilde{W}=P\rtimes W$.  
As to the action of the extended affine Hecke algebra $H(\widetilde{W})$, 
we have
$$
\iota(T^x_i f)= (T^x_i)^{-1}\iota(f)\quad(i=0,1,\cdots,n-1)
\quad \iota(\omega_x f)=\omega_x\iota(f). 
\tag{4.12}
$$
Namely, the action of $H(\widetilde{W})$ is reversed by the involution, 
denoted by the same symbol, such that
$$
\iota(T_i)=(T_i)^{-1}\quad (i=0,1,\cdots,n-1),\quad \iota(\omega)=\omega. 
\tag{4.13}
$$
Note also that $\iota(T_w)=(T_{w^{-1}})^{-1}$ for each $w\in W$. 
By this involution, the Dunkl operators and the dual Dunkl operators 
defined in (2.14), (2.15) are interchanged to each other: 
$$
 \iota(Y_i)=Y^\ast_i\quad (i=1,\cdots n). 
\tag{4.14}
$$
Formula (2.22) for the dual Dunkl operators is also obtained from (2.21) by 
this dualizing procedure. 
(Note that $\iota(D_x(u))=D_x(ut^{1-n})$. )
\par
It is easily seen that, by the involution $\iota$, the rational function 
$h^I_K=h^I_K(x,y)$ of (4.2) is transformed into 
$$
\sum_{\sigma\in\frak{S}_r}
\prod_{1\le\mu\le r}\left(
\frac{(t-1)x_{i_\mu}y_{k_{\sigma(\mu)}}}{1-tx_{i_\mu}y_{k_{\sigma(\mu)}}}
\prod_{i_\mu<j\le n}
\frac{1-x_jy_{k_{\sigma(\mu)}}}{1-tx_jy_{k_{\sigma(\mu)}}}\right)
\prod_{1\le\mu<\nu\le r} 
\frac{ty_{k_{\sigma(\mu)}}-y_{k_{\sigma(\nu)}}}
{t(y_{k_{\sigma(\mu)}}-y_{k_{\sigma(\nu)}})};
\tag{4.15}
$$
namely we have 
$$
\iota(h^I_K)= t^{-rn-\binom{r}{2}+\sum_{i\in I}i}\, x_I y_K \, h^I_K\quad(|I|=|K|=r), 
\tag{4.16}
$$
where $x_I=\prod_{i\in I} x_i$, $y_K=\prod_{k\in K}y_k$. 
If we set $\ell(I)=\ell(w_I)$ with the notation of (4.7), this can be written as
$$
\iota(h^I_K)= t^{-\ell(I)-2\binom{r}{2}}\, x_I y_K \, h^I_K\quad(|I|=|K|=r). 
\tag{4.17}
$$
By using the involution $\iota$, we see that 
these $\iota(h^I_K)$ have similar 
properties as (4.4):
$$
\alignat{2}
\text{\rm (i)}&\quad T^x_i \, \iota(h^{I}_{K}) =t \,\iota(h^{I}_{K}) 
\qquad\quad(i,i+1\in I\quad  &\text{or}\quad i,i+1\not\in I),\\ 
\text{\rm (ii)}&\quad T^x_i \, \iota(h^{I}_{K})
=t\,\iota(h^{s_i(I)}_{K})+(t-1)\iota(h^{I}_{K}) 
& (i\not\in I, i+1\in I),\tag{4.18}\\
\text{\rm (iii)}&\quad T^x_i \, \iota(h^{I}_{K})= \iota(h^{s_i(I)}_{K}) 
&(i\in I, i+1\not\in I).
\endalignat
$$
The vector subspace $\iota(V_{n,r;K})\subset\K(x,y)$ is again 
an $H(\frak{S}^x_n)$-module isomorphic to the induced representation
$\operatorname{Ind}^{H(\frak{S}_n)}
(\operatorname{triv}_{H(\frak{S}_r\times\frak{S}_{n-r})})$; 
in this case the canonical generator corresponds to 
$\iota(h^{[1,r]}_K)$.  
\par\newpage
%
\section{\S5:  Action of $D_y(u)$ on $\Pi(x,y)$}
We will now describe the action of Macdonald's $q$-difference operator 
$D_y(u)$ on the kernel $\Pi(x,y)$ in terms of the Mimachi basis.  
Recall that
$$
D_y(u) \Pi(x,y)=\Pi(x,y) F(u;y,x)
\tag{5.1}
$$
with the notation of (1.12), where 
$$
F(u;y,x)=\sum_{K\subset [1,m]}(-u)^{|K|}  t^{\binom{|K|}{2}}
\prod\Sb k\in K\\ \ell\not\in K\endSb \frac{ty_k-y_\ell}{y_k-y_\ell} 
\prod\Sb 1\le i\le n\\ k\in K\endSb \frac{1-x_iy_k}{1-tx_iy_k}.
\tag{5.2}
$$
We will give an explicit development of this $F(u;y,x)$ in terms of the Mimachi 
basis $h^{I}_{K}(x,y)$.   
\par
We start with a simple formula
$$
t^m \prod_{1\le k\le m}\frac{1-xy_k}{1-txy_k}
=1+\sum_{k=1}^m \frac{t-1}{1-txy_k}
\prod\Sb 1\le \ell\le m\\ \ell\ne k\endSb\frac{ty_k-y_\ell}{y_k-y_\ell},
\tag{5.3}
$$
which can be proved by developing the left hand side into partial fractions. 
In order to generalize this formula, let us introduce the notation
$$
\align
&a_{K|L}(y)
=\prod\Sb k\in K\\ \ell\in L\backslash K\endSb\frac{ty_k-y_\ell}{y_k-y_\ell},
\quad a_K(y)=a_{K|[1,m]}(y),
\tag{5.4}\\
&b^I_K(x,y)=
\prod\Sb i\in I\\ k\in K\endSb \frac{t(1-x_iy_k)}{1-tx_iy_k}
=t^{|I||K|}\prod\Sb i\in I\\ k\in K\endSb \frac{1-x_iy_k}{1-tx_iy_k},
\endalign
$$
for $I\subset[1,n]$ and $K\subset L\subset[1,m]$. 
With this notation, one can use (5.3) to expand 
$b^{[1,n]}_{[1,m]}=b^{[1,n]}_{[1,m]}(x,y)$
as follows:
$$
\align
&b^{[1,n]}_{[1,m]}=
\prod_{k=1}^m\frac{t(1-x_1y_k)}{1-tx_1y_k}\, b^{[2,n]}_{[1,m]} 
\tag{5.5}\\
&\quad
=\left(1+\sum_{k=1}^m \frac{t-1}{1-tx_1y_k}
\prod_{\ell\ne k}\frac{ty_k-y_\ell}{y_k-y_\ell}\right)
\, b^{[2,n]}_{[1,m]}\\
&\quad=b^{[2,n]}_{[1,m]}\\
&\qquad +\sum_{k=1}^m\left(\frac{t-1}{1-tx_1y_k}
\prod_{1<i\le n}\frac{t(1-x_iy_k)}{1-tx_iy_k}
\prod_{\ell\ne k}\frac{ty_k-y_\ell}{y_k-y_\ell}\right)
\, b^{[2,n]}_{\{1\cdots \hat{k}\cdots m\}}\\
\endalign
$$
One can use this formula repeatedly to decompose 
$b^{[2,n]}_{[1,m]}$,
$b^{[2,n]}_{\{1\cdots \hat{k}\cdots m\}}$, and so on. 
By tracing this procedure, we reach the Mimachi basis.  
\Proposition{5.1}{
With the notation of (5.4), one has
$$
b^{[1,n]}_{[1,m]}(x,y)
=\sum\Sb K\subset[1,m], I\subset[1,n]\\ |K|=|I|\endSb
a_{K}(y) h^{I}_{K}(x,y).
\tag{5.6}
$$
Furthermore, 
if $J=[n-s+1,n]$ for some $s=0,1,\cdots,n$, 
one has
$$
b^{J}_{L}(x,y)
=\sum\Sb I\subset J, K\subset L\\ |I|=|K|\endSb
a_{K|L}(y) h^{I}_{K}(x,y),
\tag{5.7}
$$
for any subset $L\subset[1,m]$. 
}
\noindent 
Note that formula (5.6) does not depend on the ordering 
of the variables $y_\ell$. 
Although it depends on the ordering of $x_j$, 
formula (5.7) is obtained from (5.6) simply by renaming 
the variables $x_1,\cdots,x_n$ as $x_{n-s+1},\cdots,x_n$
with $n$ replaced by $n-s$.  
For more details of the proof of this proposition, see \cite{MN2}. 
\par
We now try to express $F(u;y,x)$ in terms of the Mimachi basis. 
With the notation of (5.4), the function $F(u;y,x)$ can be written 
as 
$$
F(u;y,x)=\sum_{L\subset[1,m]}
(-u)^{|L|} t^{\binom{|L|}{2}-|L|n} 
a_{L}(y) b^{[1,n]}_L(x,y). 
\tag{5.8}
$$
By (5.7), we can rewrite the right-hand side into
$$
\sum\Sb I\subset[1,n],K\subset{[1,m]} \\ |I|=|K|\endSb
h^{I}_{K}(x,y) 
\sum_{K\subset L\subset[1,m]}
(-u)^{|L|} t^{\binom{|L|}{2}-|L|n} a_{L}(y) a_{K|L}(y). \tag{5.9}
$$
By the definition (5.4), it is easily seen that
$$
a_{L}(y) a_{K|L}(y)=a_{K}(y) 
a_{L\backslash K| [1,m]\backslash K}(y). 
\tag{5.10}
$$
Write $L=K\cup M$ with $M\subset[1,m]\backslash K$. 
Then the second summation of (5.9) takes the form
$$
(-u)^{|K|} t^{\binom{|K|}{2}-|K|n } a_{K}(y)
\,\sum_{M\subset[1,m]\backslash K}
(-ut^{|K|-n})^{|M|} t^{\binom{|M|}{2}}a_{M|[1,m]\backslash K}(y). 
\tag{5.11}
$$
Here we used the identity $\binom{a+b}{2}=\binom{a}{2}+\binom{b}{2}+ab$. 
This summation over $M\subset[1,m]\backslash K$ 
is nothing but the action of Macdonald's operator 
$D_{z}(ut^{|K|-n})$ 
for the variables $z=(y_\ell)_{\ell\in[1,m]\backslash K}$ 
on the constant function $1$; hence it is equal to $(ut^{|K|-n};t)_{m-|K|}$. 
Therefore, (5.11) is equal to
$$
(-u)^{|K|} (ut^{|K|-n};t)_{m-|K|} t^{\binom{|K|}{2}-|K|n } a_{K}(y).
\tag{5.12}
$$
Finally we get
$$
F(u;y,x)
=\sum_{|I|=|K|} 
(-u)^{|K|} (ut^{|K|-n};t)_{m-|K|} t^{\binom{|K|}{2}-|K|n }
a_{K}(y) h^{I}_{K}(x,y),
\tag{5.13}
$$
where the summation is taken over all pairs of indexing sets
$I\subset[1,n]$ and 
$K\subset{[1,m]}$  such that $|I|=|K|$. 
Since $(ut^{r-n};t)_{m-r}=(ut^{-n+m-1};t^{-1})_{m-r}$, 
formula (5.13) is simplified when the parameter $u$ 
is specialized to $t^{n-m+1}$. 
\Proposition{5.2}{
The function $F(u;y,x)$ is expressed in terms of the Mimachi basis as follows:
$$
\align
F(u;y,x)= \sum_{r=0}^{m\wedge n}(-u)^r(ut^{r-n};t)_{m-r} t^{\binom{r}{2}-rn}
\sum_{|I|=|K|=r} a_{K}(y) h^{I}_{K}(x,y), 
\tag{5.14}
\endalign
$$
where $I\subset[1,n]$ and $K\in[1,m]$. 
When $u=t^{n-m+1}$, this reduces to
$$
F(t^{n-m+1};y,x)=(-1)^m t^{-\binom{m}{2}} 
\sum\Sb I\subset[1,n]\\ |I|=m\endSb h^I_{[1,m]}(x,y),
\tag{5.15}
$$
if $m\le n$, and $F(t^{n-m+1};y,x)=0$ if $m>n$. 
Hence one has 
$$
D_y(t^{n-m+1})\Pi(x,y)
=(-1)^m t^{-\binom{m}{2}} \Pi(x,y) 
\sum\Sb I\subset[1,n]\\ |I|=m\endSb h^I_{[1,m]}(x,y),
\tag{5.16}
$$
if $m\le n$.
}
\noindent
We remark that formula (5.15) for $m=1$ recovers the formula
$$
1-t^n\prod_{i=1}^n\frac{1-x_iy}{1-tx_iy}=
\sum_{i=1}^n t^{n-i}\frac{1-t}{1-tx_iy} \prod_{i<j\le n}\frac{1-x_jy}{1-tx_jy}.
\tag{5.17}
$$
Note that (5.16) is also equal to $(t;t)_{n-m}^{-1} D_x(t) \Pi(x,y)$.
\par\medpagebreak\par
By the involution $\iota$ explained at the end of Section 4, 
we can easily obtain the dual version of the propositions above.  
One has only to notice that 
$$
\iota(a_{K|L}(y))=t^{-|K|(|L|-|K|)} a_{K|L}(y), \quad
\iota(b^I_K(x,y))=t^{-|I||K|} b^I_K(x,y) 
\tag{5.18}
$$
and
$$
\iota(D_y(u))=D_y(ut^{1-m}),\quad
\iota(F(u;y,x))=F(ut^{n-m+1};y,x). 
\tag{5.19}
$$
\Proposition{5.3}{
With the notation of (5.4), one has
$$
t^{-mn}b^{[1,n]}_{[1,m]}(x,y) 
=\sum\Sb I\subset[1,n], K\subset[1,m]\\ |I|=|K|\endSb 
t^{-\ell(I)-|I|(m-1)}
a_{K}(y) \, x_I y_K h^{I}_{K}(x,y). 
\tag{5.20}
$$
Furthermore, 
if $J=[n-s+1,n]$ for some $s=0,1,\cdots,n$, 
one has
$$
t^{-|J||L|}b^{J}_{L}(x,y)
=\sum\Sb I\subset J, K\subset L\\ |K|=|I|\endSb
t^{-\ell(I)-|I|(|L|-1)}
a_{K|L}(y) \, x_I y_K\,h^{I}_{K}(x,y),
\tag{5.21}
$$
for any subset $L\subset[1,m]$. 
}
\noindent
Here we used the notation 
$\ell(I)=\ell(w_I)=|I|n -\binom{|I|}{2}-\sum_{i\in I}i$. 
\Proposition{5.4}{
The function $F(u;y,x)$ is expressed in terms of the Mimachi basis as follows:
$$
\align
F(u;y,x)= \sum_{r=0}^{m\wedge n}(-u)^r(u;t)_{m-r}
\sum_{|I|=|K|=r} t^{-\ell(I)-\binom{r}{2}}a_{K}(y) 
\, x_I y_K\, h^{I}_{K}(x,y),
\tag{5.22}
\endalign
$$
where $I\subset[1,m]$ and $K\subset[1,m]$. 
When $u=1$, this reduces to
$$
F(1;y,x)=(-1)^m t^{-\binom{m}{2}} y_1\cdots y_m
\sum\Sb I\subset[1,n]\\ |I|=m\endSb t^{-\ell(I)}\,x_I \,h^I_{[1,m]}(x,y),
\tag{5.23}
$$
if $m\le n$, and $F(1;y,x)=0$ if $m>n$. 
Hence one has 
$$
D_y(1)\Pi(x,y)
=(-1)^m t^{-\binom{m}{2}} y_1\cdots y_m \Pi(x,y) 
\sum\Sb I\subset[1,n]\\ |I|=m\endSb t^{-\ell(I)}\, x_I \, h^I_{[1,m]}(x,y), 
\tag{5.24}
$$
if $m\le n$. 
}
\noindent
We remark that formula (5.23) for $m=1$ recovers the formula
$$
1-\prod_{i=1}^n\frac{1-x_iy}{1-tx_iy}=
\sum_{i=1}^n \frac{(1-t)x_iy }{1-tx_iy} \prod_{i<j\le n}\frac{1-x_jy}{1-tx_jy},
\tag{5.25}
$$
which we used in Section 3 to prove Lemma 3.3 for $m=1$. 
Note also that (5.24) is equal to 
$(t^{m-n};t)_{n-m}^{-1} D_x(t^{m-n}) \Pi(x,y)$.
%
\section{\S6:  Computation of the Dunkl operators acting on $\Pi(x,y)$}
In this section we will compute the action of the Dunkl operators and the operators 
$$
\align
&B^x_m =\sum_{1\le k_1<\cdots<k_m\le n}
 x_{k_1}\cdots x_{k_r}
(1-t^m Y_{k_1})\cdots(1-t Y_{k_m})\quad\text{and}\tag{6.1}\\
&A^x_m=\sum_{1\le k_1<\cdots<k_m\le n}
\frac{1}{ x_{k_1}\cdots x_{k_r}}
(1-Y^\ast_{k_1})\cdots(1-t^{m-1} Y^\ast_{k_m})
\endalign
$$
on $\Pi(x,y)$ by means of the Mimachi basis, to complete the proof of 
Lemma 3.3. 
{}From now on, we assume that $m\le n$. 
\par\medpagebreak\par
Recall first that the Dunkl operators satisfy the following commutation 
relations: 
$$
\overline{T}_{i} Y_{i+1} \overline{T}_i = Y_{i}, 
\quad\overline{T_i} Y_j= Y_j\overline{T_i} \quad\quad(j\ne i,i+1),
\tag{6.2}
$$
for $i=1,\cdots,n-1$. 
If $\varphi(x)$ is a function symmetric in $(x_i,x_{i+1})$,
i.e., $s_i\varphi(x)=\varphi(x)$, then 
$T_i\varphi(x)=t\varphi(x)$ and $\overline{T}_i\varphi(x)=\varphi(x)$;
hence we have  
$$
\overline{T_i}Y_{i+1}\varphi(x)=Y_i\varphi(x),\quad
\overline{T_i}Y_j\varphi(x)=Y_j \varphi(x)\quad(j\ne i,i+1). 
\tag{6.3}
$$
In the following, 
we will use this property of the Dunkl operators extensively. 
\par
For the computation of the action of the operators $B^x_m$ and $A^x_m$,
let us introduce some abbreviated notations. 
For each subset $I=\{i_1<i_2<\cdots<i_r\}$ of $[1,n]$,
we set
$$
x_I=x_{i_1}x_{i_2}\cdots x_{i_r}, \ \ 
\tau_I=\tau_{i_1}\tau_{i_2}\cdots\tau_{i_r}
\quad\text{and}\quad
Y_I=Y_{i_1} Y_{i_2}\cdots Y_{i_r}.\tag{6.4}
$$
Note that the ordering of $Y_i$'s does not matter since 
they are mutually commutative. 
Let $J=\{j_1<j_2<\cdots<j_m\}$ be a subset of $[1,n]$
and ${\bb{u}}=(u_1,u_2,\cdots,u_m)$ an $m$-tuple of parameters.  
We need to consider operators in ${\Cal D}_{q,x}[W]$ of the 
form
$$
\align
&(1-u_1 Y_{j_1})(1-u_2 Y_{j_2})\cdots(1-u_m Y_{j_m})\tag{6.5}\\
&\quad=
\sum_{r=0}^m (-1)^r 
\sum_{1\le \nu_1<\cdots<\nu_r\le m}
u_{\nu_1}\cdots u_{\nu_r} Y_{j_{\nu_1}}\cdots Y_{j_{\nu_r}}.
\endalign
$$
Such an expression will be abbreviated as follows: 
$$
(1-\bb{uY})_J 
=\sum_{I\subset J} (-\bb{u})_{I|J} Y_I,
\quad\text{with}\quad
(-\bb{u})_{I|J}=\prod\Sb 1\le \nu\le m\\ j_\nu\in I\endSb (-u_{\nu}).
\tag{6.6}
$$
With this notation, one can write formulas like 
$$
B^x_m=\sum\Sb J\subset[1,n]\\ |J|=m\endSb x_J(1-\bb{uY})_J
\quad\text{with}\quad \bb{u}=(t^m,t^{m-1},\cdots,t). 
\tag{6.7}
$$
\par
Since $\Pi(x,y)$ is a symmetric function in $x$ and $y$,
let us start with considering the action of Dunkl operators on 
{\it symmetric functions} in $x=(x_1,\cdots,x_n)$.  
Note that to consider the action of an operator  $P\in{\Cal D}_{q,x}[W]$ 
on the general symmetric function $f(x)$ is equivalent to working 
in the module ${\Cal D}_{q,x}[W]/\sum_{w\in W}{\Cal D}_{q,x}[W](w-1)$
by identifying the symbol $f(x)$ with the modulo class of $1$.  
For some time from now, $f(x)$ stands for the general symmetric 
function in $x$ in this sense. 
\Proposition{6.1}{
For the general symmetric function $f(x)$ in $x$, on has
$$
Y_I\,f(x)=\overline{T}_{w_I} \tau_{[n-r+1,n]} f(x)
\tag{6.8}
$$
for each $I\subset[1,n]$ with $|I|=r$,
where $w_I\in W$ is the permutation of  (4.6). 
Hence, for each $J\subset[1,n]$ with $|J|=m$, 
one has
$$
\align
(1-\bb{uY})_J\,f(x)&=\sum_{I\subset J} (-\bb{u})_{I|J} Y_I\,f(x) 
\tag{6.9}\\
&=\sum_{r=0}^m
\sum\Sb I\subset J\\ |I|=r\endSb (-\bb{u})_{I|J} 
\overline{T}_{w_I} \tau_{[n-r+1,n]} f(x), 
\endalign
$$
for any parameters $\bb{u}=(u_1,\cdots,u_m)$. 
}
\Proof{
By using (6.2) repeatedly, one can show that,
if $1\le i<j\le n$, then
$$
\align
&\overline{T}_i\cdots \overline{T}_{j-1} Y_j=Y_i 
\overline{T}_i^{-1}\cdots \overline{T}_{j-1}^{-1},\tag{6.10}\\
&\overline{T}_i\cdots \overline{T}_{j-1} Y_k=Y_k  
\overline{T}_i\cdots \overline{T}_{j-1}.\quad
(k<i \ \ \text{or}\ \  k>j)
\endalign
$$
Take the reduced decomposition of $w_I$ as in (4.6).  Then
it follows from (6.10) that
$$
\overline{T}_{w_I} Y_{n-r+1}\cdots Y_{n}=
Y_{i_1}\cdots Y_{i_r}\iota(\overline{T}_{w_I})
=Y_{i_1}\cdots Y_{i_r} \overline{T}_{w_I^{-1}}^{-1},\tag{6.11}
$$
namely
$$
\overline{T}_{w_I} Y_{[n-r+1,n]} \overline{T}_{w_I^{-1}}
=Y_I.
\tag{6.12}
$$
This implies 
$$
Y_I\,f(x)
=\overline{T}_{w_I} Y_{[n-r+1,n]} \overline{T}_{w_I^{-1}}f(x)
=\overline{T}_{w_I} Y_{[n-r+1,n]} f(x),
\tag{6.13}
$$
since $f(x)$ is symmetric. 
It remains to show 
$Y_{[n-r+1,n]}f(x)=\tau_{[n-r+1,n]}f(x)$. 
Since $f(x)$ is symmetric, $Y_nf(x)=\omega f(x)=\tau_n f(x)$. 
One can show inductively that 
$$
Y_{n-r+1}\cdots Y_n f(x)=\tau_{n-r+1}\cdots\tau_n f(x),
\tag{6.14}
$$
noting that the function $\tau_{n-r+1}\cdots\tau_n f(x)$ is 
symmetric in $(x_1,\cdots,x_{n-r})$ and also in 
$(x_{n-r+1},\cdots,x_n)$.
The latter half of Proposition is clear. 
}
We now apply Lemma 6.2 to the kernel $\Pi(x,y)$ regarding $y$ 
as parameters. 
Note for sure that the symbols $\tau_i$, $T_i$, $\overline{T_i}$ 
and $Y_i$ will be used below only for operators in the variables $x$. 
Since
$$
\align
\tau_{[n-r+1,n]}\Pi(x,y)
&=\Pi(x,y)\prod\Sb i\in[n-r+1,n]\\ k\in[1,m]\endSb
\frac{1-x_iy_k}{1-tx_iy_k} \tag{6.15}\\
&=\Pi(x,y) \, t^{-rm} \,  b^{[n-r+1,n]}_{[1,m]}(x,y).
\endalign
$$
we have
$$
Y_I\Pi(x,y)=\Pi(x,y) \, t^{-rm}\,  
\overline{T}_{w_I} (b^{[n-r+1,n]}_{[1,m]}(x,y)),
\tag{6.16}
$$
for each $I\subset[1,n]$ with $|I|=r$.  
To compute the action of $\overline{T}_{w_I}$ on 
$b^{[n-r+1,n]}_{[1,n]}$,
we use the expansion by the Mimachi basis:
$$
b^{[n-r+1,n]}_{[1,m]}(x,y) 
=\sum\Sb H\subset[n-r+1,n],K\subset[1,m]\\ |H|=|K|\endSb
a_{K}(y) h^H_K(x,y).
\tag{6.17}
$$
By using Theorem 4.1, one can describe the action 
$\overline{T}_{w_I}=t^{-\ell(I)} T_{w_I}$ 
on each $h^H_{K}(x,y)$.
Another notation for use: for two disjoint subsets $A,B$ of 
$[1,n]$, 
we define the {\it number of inversions} between $A$ and $B$ by
$$
\ell(A;B):=\#\{(a,b)\in A\times B \,;\, a>b\}. 
\tag{6.18}
$$
Note that, with this notation, one has 
$\ell(I)=\ell(w_I)=\ell([1,n]\backslash I; I)$. 
The action of $\overline{T}_{w_I}$ on $h^H_K(x,y)$ is then given by
$$
\overline{T}_{w_I}(h^H_K(x,y))=t^{-\ell(I)} T_{w_I}h^H_K(x,y)
=t^{-\ell([1,n]\backslash I; w_I(J))} h^{w_I(H)}_{K}(x,y).
\tag{6.19}
$$
The power of $t$ in (6.16) is obtained by subtracting the contribution
of the factors of $T_{w_I}$ which acts trivially 
(i.e., by the multiplication by $t$).  
Hence we have
$$
\align
&\overline{T}_{w_I}(b^{[n-r+1,n]}_{[1,m]}(x,y)) 
\tag{6.20}\\
&\quad=\sum\Sb H\subset[n-r+1,n], K\subset[1,m]\\ |H|=|K|\endSb
t^{-\ell([1,n]\backslash I;w_I(H))}
a_{K}(y) h^{w_(H)}_K(x,y)\\
&\quad=\sum\Sb H\subset I, K\subset[1,m]\\ |H|=|K|\endSb
t^{-\ell([1,n]\backslash I;H)}
a_{K}(y) h^{H}_K(x,y), 
\endalign
$$
where we renamed $w_I(H)$ as $H$. 
Summarizing, we get 
\Lemma{6.2}{
For each subset $I\subset[1,n]$, one has
$$
\align
&Y_I \, \Pi(x,y) 
=\Pi(x,y) \, \overline{T}_{w_I}
\prod\Sb i\in[n-r+1,n]\\ k\in[1,m]\endSb
\frac{1-x_iy_k}{1-tx_iy_k}
\tag{6.21}\\
&=\Pi(x,y)
\sum\Sb H\subset I, K\subset[1,m]\\ |H|=|K|\endSb
t^{-\ell([1,n]\backslash I; H)-|I|m}
a_{K}(y) h^H_K(x,y).
\endalign
$$
}
\par
Next we compute the action of the operator 
$(1-\bb{uY})_J=\prod_{\nu=1}^m(1-u_\nu Y_{j_\nu})$ 
on $\Pi(x,y)$ for each subset $J=\{j_1<\cdots<j_m\}$ 
of $[1,n]$  
with $|J|=m$ and for general $\bb{u}=(u_1,\cdots,u_m)$.
In view of (6.9) of Proposition 6.1, we have to compute
$$
(1-\bb{uY})_J\,\Pi(x,y)
=\Pi(x,y)\,
\sum_{r=0}^m \sum\Sb I\subset J\\ |I|=r\endSb
(-\bb{u})_{I|J} \overline{T}_{w_I}\,
\prod\Sb i\in[n-r+1,n]\\ k\in[1,m]\endSb
\frac{1-x_iy_k}{1-tx_iy_k}. 
\tag{6.22}
$$
By Lemma 6.2, the summation 
$$
\sum_{r=0}^m \sum_{I\subset J: |I|=r} 
(-\bb{u})_{I|J} \overline{T}_{w_I}\,
\prod\Sb i\in[n-r+1,n]\\ k\in[1,m]\endSb
\frac{1-x_iy_k}{1-tx_iy_k}.
\tag{6.23}
$$
is now equal to
$$
\align
&
\sum_{I\subset J} (-\bb{u})_{I|J} 
\sum\Sb H\subset I, K\subset[1,m]\\ |H|=|K|\endSb
t^{-\ell([1,n]\backslash I; H)-|I|m}
a_{K}(y) h^H_K(x,y)
\tag{6.24}\\
&=\sum\Sb H\subset J , K\subset[1,m]\\ |H|=|K|\endSb 
a_{K}(y) h^H_K(x,y) 
\sum_{H\subset I\subset J}(-\bb{u})_{I|J} 
t^{-\ell([1,n]\backslash I; H)-|I|m}
\endalign
$$
Write $I=H\cup A$ with $A\subset J\backslash H$. 
Then the summation over $I$ above can be written in the form 
$$
\align
& (-\bb{u})_{H|J} t^{-\ell([1,n]\backslash H;H)-|H|m}
\sum_{A\subset J\backslash H}
(-\bb{u})_{A|J} t^{\ell(A;H)-|A|m}
\tag{6.25}\\
&=(-\bb{u})_{H|J} t^{-\ell(H)-|H|m}
\sum_{A\subset J\backslash H}
\prod_{\nu: j_\nu\in A} (-u_\nu t^{\ell(j_\nu; H)-m})
\\
&=(-\bb{u})_{H|J} t^{-\ell(H)-|H|m}
\prod_{\nu:j_\nu\in J\backslash H} (1-u_\nu t^{\ell(j_\nu; H)-m})
\endalign
$$
If we set $u_\nu=u t^{m+1-\nu}$ ($\nu=1,\cdots,m$), 
the product at the end of (6.25) is simplified to $(u;t^{-1})_{m-|H|}$. 
This proves 
\Lemma{6.3}{
If $\bb{u}=(ut^{m},ut^{m-1},\cdots,ut)$, then 
one has
$$
\align
&\sum_{r=0}^m \sum_{I\subset J: |I|=r} 
(-\bb{u})_{I|J} \overline{T}_{w_I}\,
\prod\Sb i\in[n-r+1,n]\\ k\in[1,m]\endSb
\frac{1-x_iy_k}{1-tx_iy_k}.
\tag{6.26}\\
&
=\sum\Sb I\subset J , K\subset[1,m]\\ |I|=|K|\endSb
(-\bb{u})_{I|J}t^{-\ell(I)-|I|m}(u;t^{-1})_{m-|I|}
a_{K}(y) h^I_K(x,y),
\endalign
$$
for any subset $J\subset[1,n]$ with $|J|=m$. 
Hence
$$
\align
&(1-\bb{uY})_J\,\Pi(x,y)
\tag{6.27}\\
&=\Pi(x,y)\,
\sum\Sb I\subset J , K\subset[1,m]\\ |I|=|K|\endSb 
(-\bb{u})_{I|J}t^{-\ell(I)-|I|m}(u;t^{-1})_{m-|I|}
a_{K}(y) h^I_K(x,y). 
\endalign
$$
}
If $\bb{u}=(t^m,t^{m-1},\cdots,t)$ (i.e., $u=1$), 
the right hand side of formula (6.26)  reduces to the single term 
$$
(-1)^m t^{-\binom{m}{2}-\ell(J)} h^J_K(x,y). 
\tag{6.28}
$$
Hence by Proposition 6.1, we obtain 
\Proposition{6.4}{
Set $\bb{u}=(t^m,t^{m-1}, \cdots,t)$.  
Then for each subset 
$J=\{j_1<\cdots<j_m\}$  of $[1,n]$ with $|J|=m$, 
one has 
$$
\align
(1-\bb{uY})_J \, \Pi(x,y) 
&=(1-t^{m}Y_{j_1})(1-t^{m-1}Y_{j_2})\cdots(1-tY_{j_m}) \, \Pi(x,y)
\tag{6.29}\\
&=(-1)^m t^{-\binom{m}{2}} \Pi(x,y)\, 
t^{-\ell(J)} h^J_{[1,m]}(x,y).
\endalign
$$
Furthermore the action of $B^x_m$ on $\Pi(x,y)$ is 
expressed as follows:
$$
\align
B^x_m \, \Pi(x,y) 
&=\sum_{J\subset[1,n]:|J|=m}x_J(1-\bb{uY})_J \,\Pi(x,y) 
\tag{6.30}\\
&=(-1)^m t^{-\binom{m}{2}}\Pi(x,y) 
\sum\Sb J\subset[1,n]\\ |J|=m\endSb 
t^{-\ell(J)}x_J h^J_{[1,m]}(x,y)
\endalign
$$
}
Comparing (6.30) with the expression (5.24) of Proposition 5.4 
we obtain
\Corollary{
If $m\le n$, 
$$
B^x_m \Pi(x,y)=\frac{1}{y_1\cdots y_m} D_y(1)\,\Pi(x,y). 
\tag{6.31}
$$
}
\noindent 
This gives the proof of formula (3.17) of Lemma 3.3. 
\Remark{}{
As to the case $|J|=s< m$, the same argument as above shows
$$
\align
&(1-t^m Y_{j_1})\cdots(1-t^{m-s+1} Y_{j_s}) \Pi(x,y) 
\tag{6.32}\\
&=(-1)^s t^{-\binom{s}{2}-\ell(J)} \, 
\Pi(x,y) \sum\Sb K\subset[1,m]\\ |K|=s\endSb
a_{K}(y) h^J_K(x,y).
\endalign
$$
}
\par\medpagebreak\par
The computation of the action of $A^x_m$ can be carried out 
similarly.  
In this case, we have 
$$
Y^\ast_I f(x)=\iota(\overline{T}_{w_I})\,\tau_{[n-r+1,n]} f(x)
\quad(I\subset[1,n],\ |I|=r)
\tag{6.33}
$$
with the involution $\iota(\overline{T}_w)=(\overline{T}_{w^{-1}})^{-1}$
explained in Section 4,  
and 
$$
\align
(1-\bb{vY}^\ast)_J f(x)
&=(1-v_1Y^\ast_{j_1})(1-v_2Y^\ast_{j_2})(1-v_mY^\ast_{j_m}) f(x)
\tag{6.34}\\
&=\sum_{r=0}^m
\sum\Sb I\subset J\\ |I|=r\endSb (-\bb{v})_{I|J} 
\iota(\overline{T}_{w_I}) \, \tau_{[n-r+1,n]} f(x),
\endalign
$$
for any symmetric function $f(x)$ in $x$. 
Hence we have 
$$
\align
&(1-\bb{vY}^\ast)_J \Pi(x,y)
\tag{6.35}\\
&=\Pi(x,y) \sum_{r=0}^m
\sum\Sb I\subset J\\ |I|=r\endSb (-\bb{v})_{I|J} \iota(\overline{T}_{w_I})
\prod\Sb i\in[n-r+1,n]\\ k\in[1,m]\endSb \frac{1-x_iy_k}{1-tx_iy_k}. 
\endalign
$$
To compute this formula, we have only to dualize Lemma 6.3 
by the involution $\iota$.  
{}From Lemma 6.3, we see that, if $\bb{u}=(ut^{-m},\cdots,ut^{-1})$,
we have
$$
\align
& \sum_{r=0}^m 
\sum_{I\subset J : |I|=r } (-\bb{u})_{I|J} t^{|I|m}\iota(\overline{T}_{w_I})
\prod\Sb i\in[n-r+1,n]\\ k\in[1,m]\endSb \frac{1-x_iy_k}{1-tx_iy_k}
\tag{6.36}\\
&=\sum\Sb I\subset J, K\subset[1,m]\\ |I|=|K|\endSb 
(-\bb{u})_{I|J}(u;t)_{m-|I|} t^{|I|(m-|I|+1)}
a_{K}(y) x_I y_K h^I_K(x,y).
\endalign
$$
This implies that 
$$
\align
&(1-t^m\bb{uY}^\ast)_J \, \Pi(x,y)
=(1-uY^\ast_{j_1})\cdots(1-ut^{m-1} Y^\ast_{j_m}) \Pi(x,y)
\tag{6.37}\\
&=\Pi(x,y)\sum\Sb I\subset J, K\subset[1,m]\\ |I|=|K|\endSb
(-\bb{u})_{I|J}(u;t)_{m-|I|} t^{|I|(m-|I|+1)}
a_{K}(y) x_I y_K h^I_K(x,y).
\endalign
$$
If we set $u=1$, 
this formula reduces to
$$
(1-t^m\bb{uY}^\ast)_J \Pi(x,y)
=(-1)^m t^{-\binom{m}{2}} y_1\cdots y_m \Pi(x,y) \, 
x_J h^J_{[1,m]}(x,y).
\tag{6.38}
$$
Hence we have 
$$
\align
&A^x_m \Pi(x,y)=\sum\Sb J\subset[1,n]\\ |J|=m\endSb 
\frac{1}{x_J}(1-t^m\bb{uY}^\ast)_J\,\Pi(x,y) \tag{6.39}\\
&=(-1)^m t^{-\binom{m}{2}} y_1\cdots y_m \Pi(x,y) 
\sum\Sb J\subset [1,n]\\ |J|=m\endSb  h^J_{[1,m]}(x,y) \\
&=y_1\cdots y_m\,\Pi(x,y)\, F(t^{n-m+1};y,x)\quad(\text{by Proposition 5.2})\\
&=y_1\cdots y_m \, D_y(t^{n-m+1})\Pi(x,y). 
\endalign
$$
This completes the proof of Lemma 3.3.
\par\newpage
%
\section{\S7:  $q$-Difference raising operators}
Let $P$ be an operator in ${\Cal D}_{q,x}[W]$, and express it in the form
$$
P=\sum_{w\in W} P_w w\quad\text{with}\quad P_w\in {\Cal D}_{q,x} \quad (w\in W).  
\tag{7.1}
$$
We define the $q$-difference operator $P_{\text{sym}}$ 
by 
$$
P_{\text{sym}}=\sum_{w\in W} P_w\in {\Cal D}_{q,x}. 
\tag{7.2}
$$
Then, for any symmetric function $f(x)$, the operator $P$ acts as the 
$q$-difference operator $P_{\text{sym}}$; namely, $Pf(x)=P_{\text{sym}}f(x)$.  
By the method we used in the previous section, we can determine 
the $q$-difference operators on symmetric functions which 
arise from our raising and lowering operators.
\Theorem{7.1}{
The raising operator $B_m$ and the lowering operator $A_m$ preserve 
the ring $\K[x]^W$ of symmetric polynomials, for $m=0,1,\cdots,n$. 
Furthermore, 
these operators act on symmetric functions 
as the following $q$-difference operators : 
$$
\align
(B^x_m)_{\text{\rm sym}}
&
=\sum_{r=0}^m (-1)^r t^{\binom{r}{2}+(m-n+1)r}
\sum\Sb I\subset[1,n]\\ |I|=r\endSb 
x_I e_{m-r}(x_{[1,n]\backslash I}) a_I(x) \tau_I, 
\tag{7.3}\\
(A^x_m)_{\text{\rm sym}}
&=\sum_{r=0}^m (-1)^r t^{\binom{r}{2}}
\sum\Sb I\subset[1,n]\\ |I|=r\endSb 
x^{-1}_I e_{m-r}(x^{-1}_{[1,n]\backslash I}) a_I(x) \tau_I. 
\tag{7.4}\\
\endalign
$$
Here $e_{m-r}(x_{[1,n]\backslash I})$ and
$e_{m-r}(x^{-1}_{[1,n]\backslash I})$ 
are the elementary symmetric functions of degree $m-r$ in 
the $n-r$ variables $(x_j)_{j\in[1,n]\backslash I}$
and $(x_j^{-1})_{j\in[1,n]\backslash I}$, respectively, and
$a_I(x)=\prod_{i\in I; j\in[1,n]\backslash I}(tx_i-x_j)/(x_i-x_j)$. 
}
Before the proof of Theorem 7.1, 
we will prove that the $q$-difference operators 
$(B^x_m)_{\text{sym}}$ and $(A^x_m)_{\text{sym}}$ 
are $W$-invariant.  
Then the explicit formulas mentioned above will be determined only 
by the definition (6.1) and this $W$-invariance. 
In order to deal with $B_m$ and $A_m$ simultaneously, we consider the 
operator 
$$
\CP_m(u)=\sum_{1\le j_1<\cdots<j_m\le n} x_{j_1}\cdots x_{j_m} 
(1-ut^{m}Y_{j_1})\cdots(1-u t Y_{j_m}). 
\tag{7.5}
$$
depending on a parameter $u$.  
With the notation of Section 6, we can also write  
$$
\CP_m(u)=\sum\Sb J\subset[1,n]\\ |J|=m\endSb x_J (1-\bb{uY})_J
\quad\text{with}\quad \bb{u}=(ut^{m},ut^{m-1},\cdots,ut). 
\tag{7.6}
$$
Note that the operators $B_m$ and $A_m$ are recovered from $\CP_m(u)$ by 
$$
B_m=\CP_m(1)\quad \text{and}\quad A_m=\iota(\CP_m(t^{-m})),
\tag{7.7}
$$
where $\iota$ is the involution defined in Section 4.
Then Theorem 7.1 is an immediate consequence of the following 
Propositions 7.2 and 7.3.
\Proposition{7.2}{
For each $m=0,1,\cdots,n$, the $q$-difference operator $\CP_m(u)_{\text{\rm sym}}$
is $W$-invariant. 
Hence it preserves the ring $\K[x]^W$ of symmetric polynomials. 
}
\Proof{
We work in the quotient module 
${\Cal D}_{q,x}[W]/\sum_{w\in W}{\Cal D}_{q,x}[W](w-1)$ 
denoting by $f(x)$ the modulo class of $1$, as in Section 6.  
Then, for an operator $P\in{\Cal D}_{q,x}[W]$, 
the $W$-invariance of $(P)_{\text{\rm sym}}$ is equivalent 
to 
$$
s_i P f(x)=f(x)\quad (i=1,\cdots,n-1). 
\tag{7.8}
$$
By the definition (2.11) of $T_i$, it is also equivalent to the condition
$$
T_i P f(x)=t f(x), \quad\text{i.e.,}\quad \oT_i P f(x)=f(x),
\tag{7.9}
$$
for $i=1,\cdots,n-1$. 
We now consider the operator $\CP_m(u)$ of (7.5).
Note here that
the multiplication operators $x_1,\cdots, x_n$ have 
commutation relations
$$
\oT_i x_i \oT_i = t^{-1} x_{i+1},\quad \oT_i x_j=x_j\oT_i\quad(j\ne i,i+1)
\tag{7.10}
$$
with $\oT_i$ ($i=1,\cdots,n-1$), 
similar to (6.2) of the Dunkl operators. 
For a fixed $i$, the subsets $J$ of $[1,n]$ 
with $|J|=m$ are classified into the three groups under the action of $s_i$:
$$
\align
\text{(i)}&\quad i, i+1\not\in J,\tag{7.11}\\
\text{(ii)}& \quad i\in J,i+1\not\in J \quad \text{or}\quad i\not\in J, i+1\in J,\\
\text{(iii)}& \quad i,i+1\in J.
\endalign
$$
If $J$ satisfies (i), it is clear that 
$\oT_i x_J(1-\bb{uY})_J=x_J(1-\bb{uY})_J\oT_i$. 
The other two cases are reduced essentially to
$$
\align
&(\oT_i-1) [x_i(1-uY_i)+x_{i+1}(1-uY_{i+1})]f(x)=0, \quad\text{and} \tag{7.12}\\
&(\oT_i-1) x_i x_{i+1} (1-utY_i)(1-uY_{i+1})f(x)=0,
\endalign
$$
respectively, which can be checked directly by using (7.10), (6.2) and the 
quadratic equation
$(\oT_i-1)(\oT_i+t^{-1})=0$. 
(Note that $x_i+x_{i+1}, x_ix_{i+1}$ and $tY_i + Y_{i+1}, tY_iY_{i+1}$ commute 
with $\oT_i$. )
}
\Proposition{7.3}{
For each $m=0,1,\cdots,n$, we have
$$
\CP_m(u)_{\text{sym}}=\sum_{r=0}^m (-u)^r t^{\binom{r}{2}+(m-n+1)r} 
\sum\Sb I\subset[1,n]\\ |I|=r\endSb
x_Ie_{m-r}(x_{[1,n]\backslash I})a_I(x) \tau_I.  
\tag{7.13}
$$
}
\Proof{
By Proposition 6.1, we already have the formula
$$
\CP(u)f(x)=\sum_{r=0}^m 
\sum\Sb I\subset J\subset[1,n]\\ |I|=r, |J|=m\endSb
(-\bb{u})_{I|J}\, x_J \oT_{w_I}\tau_{[n-r+1,n]}f(x)
\tag{7.14}
$$
for the general symmetric function $f(x)$. 
Hence we see that the $q$-difference operator 
$\CP_m(u)_{\text{sym}}$ can be written in 
the form
$$
\CP_m(u)_{\text{sym}} =\sum_{r=0}^{m}\sum\Sb I\subset[1,n]\\ |I|=r\endSb 
p_I(u;x) \tau_I. 
\tag{7.15}
$$
Since $P(u)_{\text{sym}}$ is $W$-invariant by Proposition 7.2, 
we conclude that, for each $I\subset[1,n]$ with $|I|=r$, 
$p_I(u;x)=w(p_{[1,r]}(u;x))$ if $w\in W$ and $w([1,r])=I$. 
Hence we have only to determine the coefficients $p_{[1,r]}(u;x)$ for 
each $r=0,1,\cdots,m$.  
In formula (7.14), we look at the product
$$
\align
&\oT_{w_{I}} \tau_{[n-r+1,n]}f(x)\tag{7.16}\\
&=t^{-\ell(I)} (t-c(\alpha_{i_r})+c(\alpha_{i_r})s_{i_r})
\cdots(t-c(\alpha_{n-1})+c(\alpha_{n-1})s_{n-1})
\\
&\ \ 
\cdots
(t-c(\alpha_{i_1})+c(\alpha_{i_1})s_{i_1})
\cdots
(t-c(\alpha_{n-r})+c(\alpha_{n-r})s_{n-r})
\tau_{[n-r+1,n]}f(x),
\endalign
$$
for each $I=\{i_1<\cdots<i_r\}$.  
Here we used the notation 
$$
c(\alpha)=c(\ep_i-\ep_j)=\frac{1-tx_i/x_j}{1-x_i/x_j}
\tag{7.17}
$$
for each positive root $\alpha=\ep_i-\ep_j$ ($1\le i<j\le n$).
{}From expression (7.14), 
it is clear that the $q$-shift operator $\tau_{[1,r]}$ appears 
only when $I=[1,r]$ and all the terms containing $s_i$ are 
picked up in the expansion. 
In the case of $I=[1,r]$, the product of terms containing $s_i$ is 
given by
$$
c(\alpha_{r})s_{r}\cdots
c(\alpha_{n-1})s_{n-1}\cdots
c(\alpha_{1})s_{1}\cdots
c(\alpha_{n-r})s_{n-r} 
=\left(\prod_{k=1}^{N}c(\beta_k)\right) w_I,
\tag{7.18}
$$
where $N=\ell(w_{[1,r]})=r(n-r)$ and 
$\{\beta_1,\cdots,\beta_{N}\}$ is the sequence of 
positive roots associated with the reduced decomposition (4.6) of 
$w_{[1,r]}=s_r\cdots s_{n-1}\cdots s_1\cdots s_{n-r}$:
$$
\align
\{\beta_1,\cdots,\beta_{N}\} 
&=\{ \alpha_r, s_r(\alpha_{r+1}), \ldots,s_r\cdots s_{n-2}(\alpha_{n-1}),\ldots,
\tag{7.19}
\\
&\qquad\ \ \ldots,
s_r\cdots s_{n-1}\cdots s_1\cdots s_{n-r-1}(\alpha_{n-r})\}\\
&= \{\ep_i-\ep_j\,;\, 1\le i\le r, r+1\le j\le n\}. 
\endalign
$$
Namely the coefficient of $\tau_{[1,r]}$ 
arising from $\oT_{w_{[1,r]}}\tau_{[n-r+1,n]}f(x)$
is given by
$$
t^{-r(n-r)}\prod\Sb i\in[1,r]\\ j\in[r+1,n]\endSb
\frac{1-tx_i/x_j}{1-x_i/x_j}
=t^{-r(n-r)} a_{[1,r]}(x). 
\tag{7.20}
$$
In order to get the coefficient $p_{[1,r]}(u;x)$, 
we have to take the summation 
$$
\align
&p_{[1,r]}(u;x)=\sum\Sb J: [1,r]\subset J\\ |J|=m\endSb
(-\bb{u})_{[1,r]|J} x_J
t^{-r(n-r)} a_{[1,r]}(x)\tag{7.21}\\
&\quad=(-u)^r t^{mr-\binom{r}{2}-r(n-r)} 
x_1\cdots x_re_{m-r}(x_{r+1},\cdots,x_n) a_{[1,r]}(x). 
\endalign
$$
Note here that 
$(-\bb{u})_{[1,r]|J}=(-u)^r t^{mr-\binom{r}{2}}$ for any $J$ with 
$[1,r]\subset J$, $|J|=m$. 
{}From the $W$-invariance of $\CP_m(u)_{\text{sym}}$, we have 
$$
p_{I}(u;x)=(-u)^r t^{mr-\binom{r}{2}-r(n-r)} 
x_Ie_{m-r}(x_{[1,n]\backslash I}) a_{I}(x),
\tag{7.22}
$$
for each $I\subset[1,n]$ with $|I|=r$. 
This proves (7.13).
}
Let us see some consequences of Proposition 7.3.  
We apply the operator $\CP_m(u)$ to the constant function $1$.  
By the definition (7.5) of $\CP_m(u)$, we have 
$$
\CP_m(u)(1) = (ut;t)_m e_m(x). 
\tag{7.23}
$$
On the other hand, by (7.13) we have
$$
\CP_m(u)(1)=\sum_{r=0}^m (-u)^r t^{\binom{r}{2}+(m-n+1)r} 
\sum\Sb I\subset[1,n]\\ |I|=r\endSb
x_Ie_{m-r}(x_{[1,n]\backslash I})a_I(x).
\tag{7.24}
$$
By comparing the coefficients of $u^r$ in (7.23) and (7.24), we obtain
\Corollary{
If $0\le r\le m\le n$, we have
$$
\align
\sum\Sb I\subset[1,n]\\ |I|=r\endSb
x_Ie_{m-r}(x_{[1,n]\backslash I})
\prod\Sb i\in I \\ j\not\in I\endSb \frac{tx_i-x_j}{x_i-x_j}
&=t^{(n-m)r}\left[ \matrix m\\ r \endmatrix\right]_t e_m(x),\tag{7.25}\\
\sum\Sb I\subset[1,n]\\ |I|=r\endSb
x_Ie_{m-r}(x_{[1,n]\backslash I})
\prod\Sb i\in I \\ j\not\in I\endSb \frac{x_i-tx_j}{x_i-x_j}
&=\left[ \matrix m\\ r \endmatrix\right]_t e_m(x).\tag{7.26}
\endalign
$$
}
\noindent
Formula (7.26) is obtained from (7.25) by the transformation $t\to t^{-1}$.
These formulas give a refinement of the equality 
$$
B_m(1)=(t;t)_m e_m(x)=J_{(1^m)}(x) \tag{7.27}
$$ 
We also remark that formulas (2.21) and (2.22) in Section 2 can be obtained
from Proposition 7.3 for $m=n$. 
In fact, we have 
$$
\align
P_n(u)&=x_1\cdots x_n (1-ut^n Y_1)\cdots(1-utY_n)\quad \text{and}
\tag{7.28} \\
P_n(u)_{\text{\rm sym}}&=x_1\cdots x_n 
\sum_{r=0}^n (-ut)^r t^\binom{r}{2}
\sum\Sb I\subset[1,n]\\ |I|=r\endSb 
a_I(x) \tau_I. 
\endalign
$$ 
This implies (2.21). 
\par\medpagebreak
\par\newpage
\section{\S8:  A double analogue of the multinomial coefficients}
In this section, we will give an application of our results to combinatorics 
of the multinomial coefficients, in terms of the so-called 
modified Macdonald polynomials. 
The modified Macdonald polynomials were introduced by A.M.\,Garsia 
and M.\,Haiman \cite{GH} in their study of a graded representation 
model for Macdonald polynomials (see also \cite{Ma1}, p.358). 
\par\medpagebreak
Following \cite{GH},  we define the {\it modified Macdonald polynomials in
infinite number of variables} by
$$
\widetilde{P}_{\ld}(x;q,t)=P_\ld(\frac{x}{1-t};q,t)\quad \text{and} \quad
\widetilde{J}_{\ld}(x;q,t)=J_\ld(\frac{x}{1-t};q,t),
\tag{8.1}
$$
in the $\lambda$-ring notation. 
Given a symmetric function $f(x)=f(x_1,x_2,\cdots)$ in infinite variables
$x=(x_1,x_2,\cdots)$, 
the symbol $f(\frac{x}{1-t})$ in the $\lambda$-ring notation
stands for the symmetric function $f(y(x))$ obtained by  
the transformation of variables $y(x)=(x_it^j)_{i\ge 1, j\ge 0}$. 
In infinite variables, the symmetric function 
$f(x)$ can be written uniquely in the form $f(x)=\varphi(p_1(x),p_2(x),\cdots)$ 
as a polynomial of the power sums  
$p_k(x)=\sum_{j=1}^\infty x_j^k$ ($k=1,2,\cdots$).  
Then the symbol $f(\frac{x}{1-t})$ represents the symmetric function
$$
f(\frac{x}{1-t})=\varphi(\frac{p_1(x)}{1-t}, \frac{p_2(x)}{1-t^2},\cdots). 
\tag{8.2}
$$
obtained by the transformation $p_k(x) \to p_k(x)/(1-t^k)$
($k=1,2,\cdots$).
When we consider the modified Macdonald polynomials in $n$ variables 
$x=(x_1,\cdots,x_n)$, each of $\widetilde{P}_\ld(x;q,t)$ and 
$\widetilde{J}_\ld(x;q,t)$ should be understood as the one obtained from 
the corresponding symmetric function in infinite variables by setting 
$x_{n+1}=x_{n+2}=\cdots=0$. 
It follows from the orthogonality property of Macdonald polynomials
(see (1.7) or \cite{Ma1}, (VI.4.13)) that
$$
\sum_\ld b_\ld \widetilde{P}_\ld(x)\, P_\ld(y)=
\widetilde{\Pi}(x,y), \ \ \text{where}\ \ 
\widetilde{\Pi}(x,y)=\prod_{i,j}
\frac{1}{(x_i y_j;q)_\infty}.
\tag{8.3}
$$
\par
An advantage of modified Macdonald polynomials is that they have 
nice transition coefficients with classical Schur functions $s_\ld(x)$:
$$
\widetilde{J}_\mu(x;q,t)=\sum_{\ld} K_{\ld,\mu}(q,t) \, s_\ld(x),
\tag{8.4}
$$
where $K_{\ld,\mu}(q,t)$ are the double Kostka coefficients.  
By Theorem 3.2, we already know that $K_{\ld,\mu}(q,t)\in\Z[q,t]$ 
for all $\ld$ and $\mu$. 
\par
We now introduce a family of functions $B_{\ld,\mu}(q,t)$ via 
decomposition
$$
\widetilde{J}_\ld(x;q,t)=\sum_{\ld} B_{\ld,\mu}(q,t) \, m_\mu(x) 
\tag{8.5}
$$
in terms of monomial symmetric functions.  
Note that $B_{\ld,\mu}(q,t)=0$ unless $|\ld|=|\mu|$. 
Using the Kostka numbers $K_{\ld,\mu}$ defined by 
$$
s_\ld(x)=\sum_{\mu} K_{\ld,\mu}\, m_\mu(x),
\tag{8.6}
$$
one can express the coefficient $B_{\ld,\mu}(q,t)$ as 
$$
B_{\ld,\mu}(q,t)=\sum_{\nu} K_{\nu,\ld}(q,t) \, K_{\nu,\mu}.
\tag{8.7}
$$
\Theorem{8.1}{
For any partitions $\ld$ and $\mu$ of a given natural number $n$, we have
\roster
\item $B_{\ld,\mu}(q,t)\in \Z[q,t]$,
\item $B_{\ld,\mu}(1,1)={\displaystyle\frac{n!}{\mu_1! \mu_2! \cdots}}$, 
\item 
$B_{(n),\mu}(q,t)=q^{n(\mu')}
{\displaystyle\frac{(q;q)_n}{(q;q)_{\mu_1}(q;q)_{\mu_2}\cdots}}$,
\item
$B_{\ld',\mu}(q,t)=q^{n(\ld')}t^{n(\ld)} B_{\ld,\mu}(t^{-1},q^{-1})$
\quad(Duality),
\endroster
where $n(\ld')=\sum_{s\in\ld} a(s)$ and 
$n(\ld)=\sum_{s\in\ld} \ell(s)$. 
}
\proclaim{Conjecture 8.2?}
$B_{\ld,\mu}(q,t)\in \N[q,t]$ for any partitions $\ld$ and $\mu$.
\endproclaim
\noindent
We remark that Conjecture 8.2? follows from the positivity conjecture 
of Macdonald \cite{Ma1}, (VI.8.18?) on the double Kostka polynomials. 
Hence it would be natural to consider the polynomials $B_{\ld,\mu}(q,t)$ 
as a two-parameter deformation of the multinomial coefficients. 
\Proofof{Theorem 8.1}{
The first statement follows from Theorem 3.2.(2) since 
$K_{\nu,\mu}\in\N$.
As for the second statement, let us remark that 
$K_{\nu,\ld}(1,1)=f^\nu$ is the number of standard Young tableaux of 
shape $\nu$ (\cite{Ma1}, (VI.8.16)).
Hence 
$$
B_{\ld,\mu}(1,1)=\sum_{\nu}f^\nu  K_{\nu,\mu} =
\frac{n!}{\mu_1!\mu_2!\cdots}. \tag{8.8}
$$
The last equality is well-known and in fact follows for example from 
the Robinson-Schensted-Knuth correspondence (see e.g. \cite{S}).  
Similarly, if $\ld=(n)$, then one has $K_{\nu,(n)}(q,t)=K_{\nu,(1^n)}(q)$ 
(see \cite{Ma1}, p.362). 
Hence
$$
B_{(n),\mu}(q,t)=\sum_{\nu} K_{\nu,(1^n)}(q) K_{\nu,\mu}
=q^{n(\mu')} \frac{(q;q)_n}{(q;q)_{\mu_1}(q;q)_{\mu_2}\cdots}.
\tag{8.9}
$$
The last equality is also well-known and is proven for example in 
\cite{DJKMO} or \cite{K}, \S2.4. 
Finally, statement (4) follows from the corresponding duality theorem 
for the double Kostka polynomials (see \cite{Ma1}, (VI.8.15) and (VI.8.5)).
}
%

\enddocument